\documentclass{article}[12pt]
\usepackage{
latexsym,
amsthm,amsfonts,amssymb,
epsfig,psfrag,array
}
\makeatletter

\long\def\@makecaption#1#2{%
  \vskip\abovecaptionskip
  \sbox\@tempboxa{{\bf #1.} #2}%
  \ifdim \wd\@tempboxa >\hsize
    \unhbox\@tempboxa\par
  \else
    \global \@minipagefalse
    \hbox to\hsize{\hfil\box\@tempboxa\hfil}%
  \fi
  \vskip\belowcaptionskip}

\def\tcaption #1{\refstepcounter{table}
\begingroup
  \@parboxrestore
  \normalsize
  \@makecaption{\fnum@table}{
   #1}\par
\endgroup}

\makeatother

\begin{document}


\newtheorem{lemma}{Lemma}
\newtheorem{col}{Corollary}
\newtheorem{theorem}{Theorem}
\renewcommand{\proofname}{Proof}

\def \H{{\mathbb H}}
\def \N{{\mathbb N}}
\def \Sp{{\mathbb S}}
\def \E{{\mathbb E}}
\def \Z{{\mathbb Z}}

\vspace{-10pt}

\begin{center}
{\LARGE  \bf
Reflection subgroups \\
of Euclidean reflection groups.
}
\vspace{10pt}

{\Large A.~Felikson, P.~Tumarkin}
\end{center}

{\large

\section{Introduction}
Let $X$ be an $n$-dimensional Euclidean, spherical, or hyperbolic space.
A convex polytope  in $X$ is called a {\it Coxeter polytope}
if all its dihedral angles are integer submultiples of $\pi$.
A group generated by reflections with respect to the facets of
a Coxeter polytope is discrete; a fundamental domain
of this group is the initial Coxeter polytope.

From now on by polytope we mean a finite volume polytope in
the $n$-dimensional Euclidean space $\E^n$.
A group generated by  reflections with respect to the facets of
 Euclidean Coxeter polytope we call a Euclidean {\it reflection group}.
Finite groups generated by reflections  we call {\it spherical}
groups.

A classification of reflection subgroups of spherical reflection groups
may be deduced from~\cite{Dyn}.
A classification of reflection subgroups of Euclidean and hyperbolic
reflection groups is still incomplete.
Papers~\cite{Fel4-9},~\cite{Fel3},~\cite{Ruth} and~\cite{Klim}
are devoted to reflection subgroups of hyperbolic reflection
groups with simplicial fundamental domains.

In this paper, we classify reflection subgroups of discrete
Euclidean reflection groups.

\vspace{8pt}
In section~\ref{sph} we make use results of~\cite{Dyn} to classify
reflection subgroups of indecomposable spherical reflection groups.
In fact, any reflection subgroup in any decomposable reflection group
$G$ is a direct product of reflection subgroups of indecomposable
components of $G$ (Lemma~\ref{direct product}).  Hence,  it is
sufficient to describe 
subgroups of indecomposable spherical and Euclidean reflection groups.
All the indecomposable Euclidean compact Coxeter polytopes are simplices
(see~\cite{Cox}). Thus, to obtain a classification of reflection subgroups
of Euclidean reflection groups we only need to classify reflection
subgroups of groups generated by reflections in the facets of Euclidean
simplices.

In Section~\ref{simpl}, we classify all indecomposable finite index
reflection subgroups of indecomposable Euclidean reflection groups.
We also prove that any indecomposable reflection subgroup is determined
by its index up to an automorphism of the whole group.  Furthermore, in
Section~\ref{General} we give a general description of reflection
subgroups in terms of affine root systems.  In Section~\ref{Max}, we
classify decomposable maximal reflection subgroups.  In Section~\ref{Inf},
we consider infinite index reflection subgroups.

The authors are grateful to E.~B.~Vinberg for helpful discussions, and
to the anonymous referee for pointing out numerous inaccuracies in
the preliminary version of the paper.


\section{Definitions and notation}

A hyperplane $\mu$ is called a {\it mirror} of a reflection group if
the group contains a reflection in $\mu$. Mirrors of a reflection group
decompose the space into connected components ({\it fundamental
chambers}); any fundamental chamber is a fundamental domain for the
group action. Reflections in facets of any fundamental chamber
of a reflection group generate the whole group.

Let $H$ be a finite index reflection subgroup of a reflection group $G$.
Then a fundamental chamber of $H$ consists of several copies of a
fundamental chamber of $G$. If two such copies have a facet in common then
they are symmetric to each other with respect to this facet.

To describe Coxeter polytopes we use {\it Coxeter diagrams}:
nodes  $v_1$, $\dots$, $v_k$ of the diagram correspond
to the facets $f_1,...,f_k$ of the polytope;
the nodes $v_i$ and $v_j$ are joined by  an $(m_{ij}-2)$-fold
edge if the dihedral angle formed up by  $f_i$ and $f_j$ is equal to
$\frac{\pi}{m_{ij}}$ (if  $f_i$  is orthogonal to $f_j$ then
the nodes $v_i$ and $v_j$ are disjoint);
the nodes $v_i$ and $v_j$ are joined by a bold edge
if  $f_i$  is parallel to $f_j$.

Let $\Sigma$ be a Coxeter diagram with
$k$ nodes $v_1$,...,$v_k$.
Denote by $Gr_{\Sigma}=( g_{ij})$ a symmetric  $k\times k$ matrix with
$g_{ii}=1$ ($i=1,...,k$),
$g_{ij}=-\cos(\frac{\pi}{m_{ij}})$  when  $i\ne j$
and  $v_i$ is connected with $v_j$ by  an $(m_{ij}-2)$-fold
edge, $g_{ij}=-1$  when  $i\ne j$
and  $v_i$ is connected with $v_j$ by a bold edge.

A connected Coxeter diagram  $\Sigma$ is called {\it elliptic} if
$Gr_{\Sigma}$ is positively defined; a connected diagram $\Sigma$ is called
{\it parabolic} if  $Gr_{\Sigma}$ is degenerate and any subdiagram of
 $\Sigma$ is elliptic. A Coxeter polytope is called {\it indecomposable}
if its Gram matrix is indecomposable (or, similarly, its Coxeter diagram
is connected). Connected elliptic Coxeter diagrams coincide
with Coxeter diagrams of indecomposable spherical Coxeter
simplices, and connected parabolic diagrams coincide with Coxeter diagrams
of Euclidean Coxeter simplices.  A Coxeter diagram of any
compact Euclidean Coxeter polytope is a disjoint union of connected
parabolic diagrams.  Table~\ref{cox} contains the list of Coxeter diagrams
of indecomposable spherical and Euclidean Coxeter simplices.  We also use
the notation $B_1=C_1=A_1$, $D_2=2A_1$ and $D_3=A_3$ (see~\cite{Dyn}).

%
\begin{table}
\tcaption{{\bf Coxeter diagrams.} 
Connected elliptic and parabolic Coxeter diagrams are
listed in left and right columns respectively.
Special nodes are colored in white.
\\  }
\label{cox}
\begin{center}
\begin{tabular}{|cc@{\quad}|cc|}
\hline
\raisebox{0pt}{${ A_n}$ $(n\ge 1)$}  & 
\raisebox{0pt}{\epsfig{file=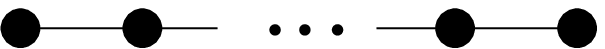,width=0.2\linewidth}}&
\multicolumn{2}{c|}{
\begin{tabular}{cc}
\raisebox{0pt}[20pt][5pt]{${ \widetilde A_1}$} & 
\raisebox{0pt}[20pt][5pt]{\epsfig{file=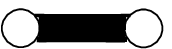,width=0.052\linewidth}}\\
\raisebox{3pt}[15pt][14pt]{${ \widetilde A_n}$ $(n\ge 2)$}  & 
\raisebox{-8pt}[25pt][7pt]{\epsfig{file=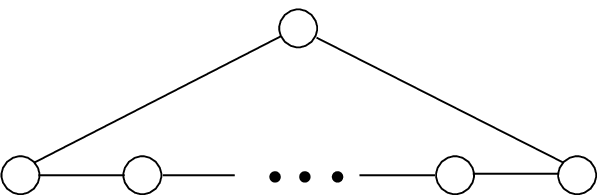,width=0.2\linewidth}}
\end{tabular}
}\\
\hline
\raisebox{-1pt}[23pt][7pt]{${ B_n=C_n}$} & 
\raisebox{-7pt}[23pt][7pt]{\epsfig{file=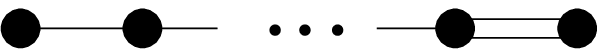,width=0.2\linewidth}}&
\raisebox{7pt}[23pt][7pt]{${ \widetilde B_n}$ $(n\ge 3)$}  & 
\raisebox{-0pt}[30pt][7pt]{\epsfig{file=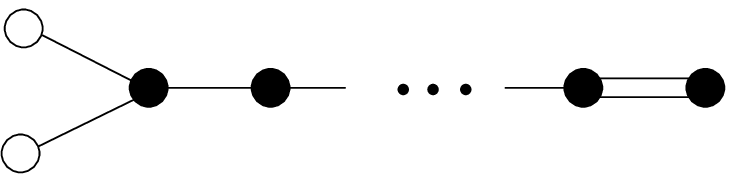,width=0.2\linewidth}}\\
\cline{3-4}
\raisebox{7pt}[15pt][7pt]{ $(n\ge 2)$} & 
&
\raisebox{-0pt}[15pt][7pt]{${ \widetilde C_n}$ $(n\ge 2)$}  & 
\raisebox{-0pt}[15pt][7pt]{\epsfig{file=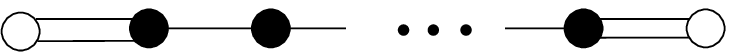,width=0.2\linewidth}}\\
\hline
\raisebox{7pt}[23pt][7pt]{${ D_n}$ $(n\ge 4)$} & 
\raisebox{0pt}[30pt][7pt]{\epsfig{file=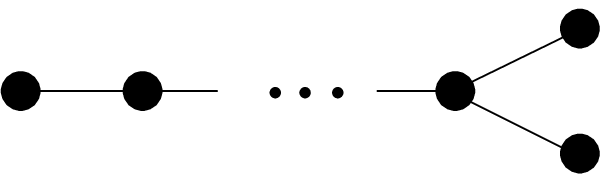,width=0.2\linewidth}}&
\raisebox{7pt}[23pt][7pt]{${ \widetilde D_n}$ $(n\ge 4)$} & 
\raisebox{0pt}[30pt][7pt]{\epsfig{file=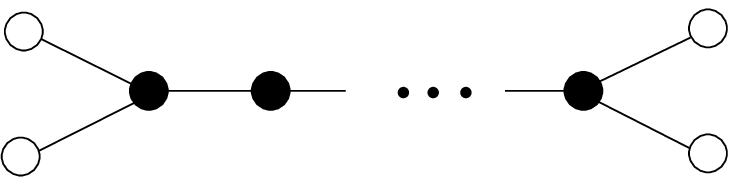,width=0.2\linewidth}}\\
\hline
\raisebox{0pt}[15pt][7pt]{${ G_2}$}  & 
\raisebox{0pt}[15pt][7pt]{\epsfig{file=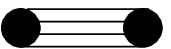,width=0.052\linewidth}}&
\raisebox{0pt}[15pt][7pt]{${ \widetilde G_2}$} & 
\raisebox{0pt}[15pt][7pt]{\epsfig{file=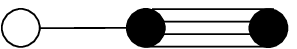,width=0.082\linewidth}}\\
\hline
\raisebox{0pt}[15pt][7pt]{${ F_4}$}  & 
\raisebox{0pt}[15pt][7pt]{\epsfig{file=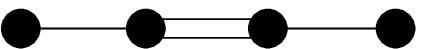,width=0.13\linewidth}}&
\raisebox{-0pt}[15pt][7pt]{${ \widetilde F_4}$} & 
\raisebox{-0pt}[15pt][7pt]{\epsfig{file=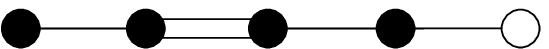,width=0.16\linewidth}}\\
\hline
\raisebox{8pt}[15pt][7pt]{${ E_6}$}  & 
\raisebox{0pt}[30pt][7pt]{\epsfig{file=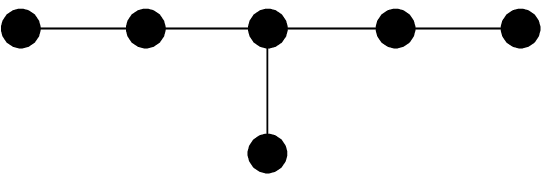,width=0.16\linewidth}}&
\raisebox{8pt}[35pt][7pt]{${ \widetilde E_6}$} & 
\raisebox{-8pt}[40pt][17pt]{\epsfig{file=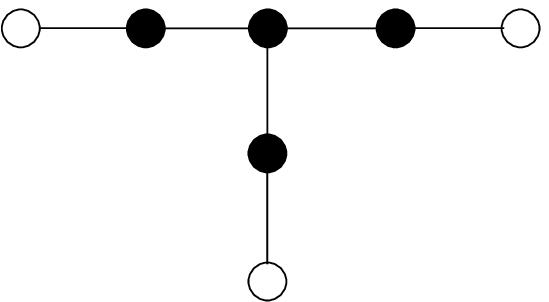,width=0.16\linewidth}}\\
\hline
\raisebox{5pt}[25pt][7pt]{${ E_7}$}  & 
\raisebox{-0pt}[30pt][7pt]{\epsfig{file=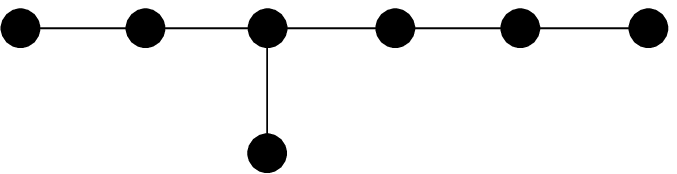,width=0.2\linewidth}}&
\raisebox{5pt}[25pt][7pt]{${ \widetilde E_7}$} & 
\raisebox{-0pt}[25pt][7pt]{\epsfig{file=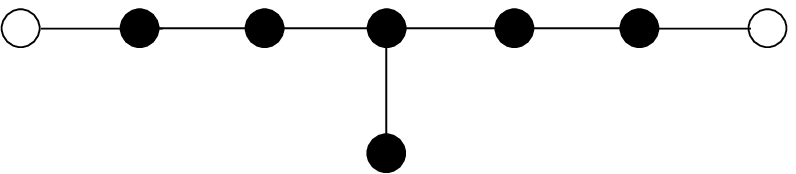,width=0.22\linewidth}}\\
\hline
\raisebox{5pt}[25pt][7pt]{${ E_8}$}  & 
\raisebox{-0pt}[30pt][7pt]{\epsfig{file=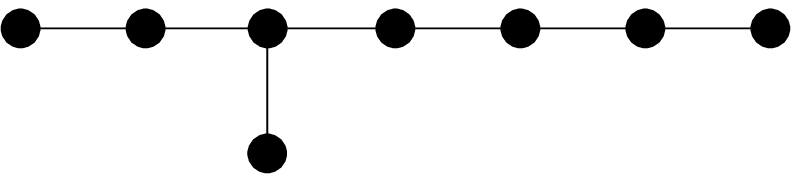,width=0.24\linewidth}}&
\raisebox{5pt}[25pt][7pt]{${ \widetilde E_8}$} & 
\raisebox{-0pt}[25pt][7pt]{\epsfig{file=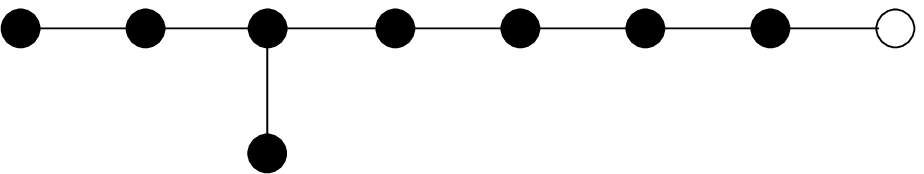,width=0.25\linewidth}}\\
\hline
\raisebox{0pt}[15pt][7pt]{${ H_3}$}  & 
\raisebox{0pt}[15pt][7pt]{\epsfig{file=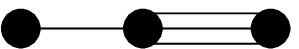,width=0.1\linewidth}}&
& 
\\
\cline{1-2}
\raisebox{0pt}[15pt][7pt]{${ H_4}$}  & 
\raisebox{0pt}[15pt][7pt]{\epsfig{file=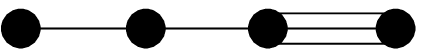,width=0.13\linewidth}}&
& 
\\
\hline

\end{tabular}
\end{center}

\end{table}


Denote by $\Sigma(P)$ the Coxeter diagram of a Coxeter polytope $P$
and by $G_P$ the group generated by reflections in the
facets of $P$. The matrix $Gr_{\Sigma(P)}$ coincides with
the Gram matrix of $P$.

A survey about reflection groups and their fundamental polytopes
may be found in~\cite{29}.

Let $P$ be a Euclidean Coxeter polytope. There exists at least one vertex
$V$ of $P$ such that the stabilizer of $V$ contains linear parts of all
elements of $G_P$ (see~\cite[Ch. 6]{Kac}). Vertices of that type are
called {\it special vertices of polytope $P$}.

Now let $P$ be an indecomposable Euclidean Coxeter polytope, i.e.
a Coxeter simplex; let $V$ be any special vertex of $P$ and $\mu$ be a
facet opposite to $V$.  Let $v$ be a node of $\Sigma(P)$ corresponding to
$\mu$.  We call $v$ a {\it special node of} $\Sigma(P)$.
Notice that $\Sigma(P)$ may contain several special vertices.  However,
the diagram $\Sigma(P)\setminus v$ does not depend on the choice of
special vertex (in other words, any two diagrams of this type are
equivalent under an automorphism of the diagram $\Sigma(P)$).  Denote
$\Sigma(P)\setminus v$ by  $\Sigma'(P)$.  Denote by $G'_P(V)$ the group
generated by reflections in the facets of $P$ different from $\mu$.
This group does not depend on the choice of special vertex $V$ (i.e.
for any pair of special vertices   $V_1$ and $V_2$ of $P$ there exists an
automorphism of $G_P$ taking $V_1$ to $V_2$.  When the choice of $V$
is not important 
we write $G'_P$ instead of $G'_P(V)$.  Note that $\Sigma'(P)$ is a Coxeter
diagram of a fundamental chamber of $G'_P$.

A subgroup $G_P\subset G_F$ is called {\it maximal} if there is no simplex
$T$ such that $G_P\subset G_T\subset G_F$.

A reflection group is called {\it indecomposable} if its fundamental
Coxeter polytope is indecomposable. Any reflection group is a direct
product of several indecomposable reflection groups. We call these factors
{\it components} of a reflection group.

\begin{lemma}
\label{direct product}
Any reflection subgroup $H$ of decomposable reflection group $G$ is a
direct product of reflection subgroups of indecomposable components of
$G$.

\end{lemma}

\begin{proof}
Let $G=G_1\times G_2\times\dots\times G_k$, and let
$\left\{r_1^i,\dots,r_{l_i}^i\right\}$ be reflections generating the group
$G_i$. Denote by $\left[G,G\right]$ the commutator subgroup of $G$.
Furthermore, notice that
$$G/\left[G,G\right]\cong\Z_2^m,$$
where $M$ is the number of conjugacy classes of reflections in $G$.
In other words, any set of reflections generating $G$ must contain
representatives of all conjugacy classes of reflections in $G$. Hence,
any reflection $r\in H$ is conjugated in $G$ to some $r_j^i\in G_i$. Since
$G_i\vartriangleleft G$, we obtain that $r\in G_i$, and the lemma is
proved.

\end{proof}

\section{Spherical reflection subgroups}
\label{sph}

In paper~\cite{Dyn} Dynkin classified regular semisimple subalgebras of
semisimple Lie algebras. For this he listed all root subsystems in finite
root systems. This problem is very close to the  classification
of reflection subgroups of finite reflection groups. However, it is not the
same.

In this section we make use the results of~\cite{Dyn} to classify
reflection subgroups of finite reflection groups.

\subsection{Root subsystems}

Let $G$ be a finite reflection group different from $H_3$, $H_4$ and
$G_2^{(m)}$ (if $m\ne 2,3,4,6$). Then $G$ may be thought as a Weyl group
of some finite root system $\Delta$ (see~\cite{Bur}). Mirrors of $G$ are
hyperplanes on which the roots of $\Delta$ vanish. A type of $G$ coincide with
the type of the root system $\Delta$.

Any reflection subgroup $H\subset G$ is a Weyl group of some root system
$\Delta_1\subset\Delta$. The root system $\Delta_1$ consists of those
roots of $\Delta$ which vanish on mirrors of $H$.
Conversely, a Weyl group of any root system $\Delta_1\subset\Delta$
is a reflection subgroup of $G$.

In this way we obtain a one-to-one correspondence between reflection
subgroups $H\subset G$ and root systems $\Delta_1\subset\Delta$.

A root system $\Delta_1\subset\Delta$ is called {\it a root subsystem}
if the following holds:  $$ {\mbox{if }}\
\alpha, \beta\in\Delta_1 \ {\mbox {and }} \ \alpha+\beta\in\Delta,\ {\mbox
{then} \ }\ \alpha+\beta\in\Delta_1 \eqno{(*)} $$

Root subsystems of finite root systems are classified in~\cite{Dyn}.
So, for each finite reflection group $G$ we only need to list
root systems $\Delta_1\subset\Delta$ which are not root subsystems.

Condition $(*)$ may be reformulated in the following way.
Let $\Delta_1\subset\Delta$ be root systems, and let
 $\{\alpha_1,\dots,\alpha_{n}\}$ be simple roots of $\Delta_1$.
The root system $\Delta_1$ is a root subsystem of $\Delta$ if and
only if the following holds: $$
\alpha_i-\alpha_j\notin\Delta \mbox{\ if } i\ne j.\eqno(**) $$

What does condition $(**)$ mean from geometric point of view?
Let $\Delta_1\subset\Delta$ be root systems, let
$\alpha_i,\alpha_j\in\Delta_1$ be simple roots of $\Delta_1$, and let
$\alpha_i-\alpha_j\in\Delta$. Consider the root system generated
by $\alpha_i$ and $\alpha_j$. It should coincide with one of $2A_1, A_2,
C_2$ or $G_2$.  For each of these systems the angle formed up by roots
$\alpha_i$ and $\alpha_j$ is cut by the line orthogonal to the root
$\alpha_i-\alpha_j$ (simple roots are the outward normal vectors to the
sides of the angle). The angle formed up by simple roots can not be
acute, so we have
$$ (\alpha_i-\alpha_j)^2=
\alpha_i^2+\alpha_j^2-2(\alpha_i,\alpha_j)\ge\alpha_i^2+\alpha_j^2.
$$
Since any finite root system contains roots of at most two
different lengths, the equality holds if and only if $\alpha_i$ and
$\alpha_j$ are short roots of $\Delta$,
$(\alpha_i,\alpha_j)=0$ when $\Delta\ne G_2$, and
$(\alpha_i,\alpha_j)=-\frac{1}{2}$ when $\Delta=G_2$.

In particular, we obtain the following
\begin{lemma}
\label{cut}
Let $\Delta_1\subset\Delta$ be root systems, let $\Pi_1$ and
$\Pi$ be simple root systems of  $\Delta_1$ and
$\Delta$ respectively. Suppose that one of the following holds:

1) all the real roots of $\Delta$ have the same length;

2) $\Pi_1\setminus\Pi$ contains no short root.

Then $\Delta_1$ is a root subsystem of $\Delta$.

\end{lemma}

\begin{col}
\label{sl}
All reflection subgroups of the reflection groups $A_n$, $D_n$, $E_6$, $E_7$
and $E_8$ are listed in~\cite{Dyn}.

\end{col}

Hence, we only need to list root systems $\Delta_1\subset\Delta$ that 
are not root subsystems, when $\Delta=B_n$ (or $C_n$), $F_4$ and $G_2$,
as well as reflection subgroups of reflection groups  $H_3$, $H_4$ and
$G_2^{(m)}$ (when $m\ge 5,\ m\ne 6$).

\subsection{Subgroups of $B_n$}
\label{b_n-sec}

Consider the group $G=B_n$ as the Weyl group of the root system
$\Delta=B_n$. By~\cite[Table~9]{Dyn} (see also~\cite{On}),
the root system $B_n$ contains subsystems of the type
$$
A_{k_1}+\dots+A_{k_s}+D_{m_1}+\dots+D_{m_r}+B_m,
$$
$$
\sum\limits_{i=1}^s(k_i+1)+\sum\limits_{i=1}^r m_i+m\le n,\quad
k_1\ge\dots\ge k_s\ge 0,\quad  m_1\ge\dots\ge m_r>1,\quad m\ge 0.
$$

From the other hand, the group $B_n$ may be considered as the Weyl group
of the root system $C_n$. Changing lengths of roots in root subsystems of
$C_n$ (see~\cite[Table~9]{Dyn}), we see that $B_n$ contains root
subsystems of the type
$$
A_{k_1}+\dots+A_{k_s}+B_{l_1}+\dots+B_{l_p},
$$
where \qquad
$
\sum\limits_{i=1}^s(k_i+1)+\sum\limits_{i=1}^p l_i\le n,\quad
k_1\ge\dots\ge k_s\ge 0,\quad  l_1\ge\dots\ge l_p>0.
$

An elementary check shows that all root systems $\Delta_1\subset\Delta$
consist of components of the systems described above. More
precisely, the following holds:

\begin{lemma}
\label{b_n}
Let $\Delta=B_n$ and $\Delta_1\subset\Delta$.
Then $\Delta_1$ is of the type
$$
A_{k_1}+\dots+A_{k_s}+B_{l_1}+\dots+B_{l_p}+D_{m_1}+\dots+D_{m_r},
$$
where\qquad
$
\sum\limits_{i=1}^s(k_i+1)+\sum\limits_{i=1}^p l_i+\sum\limits_{i=1}^r
m_i \le n,\quad k_1\ge\dots\ge k_s\ge 0,\quad l_1\ge\dots\ge l_p>0,
$
\vspace{-3pt}
\begin{center}
$
m_1\ge\dots\ge m_r>1.
$
\end{center}

For any two root systems of the same type there exists an automorphism of
$\Delta$ taking one system to another.

\end{lemma}

Notice that components of the type $A_1$ and $D_2$ consist of long
roots, and components of the type $B_1$ consist of short roots.
Furthermore, two different root systems $\Delta_1,\,\Delta_2\subset\Delta$
may have Weyl groups $H'$ and $H''$ of the same type. For example,
Weyl groups of root systems $\Delta_1=A_3+D_2$ and $\Delta_2=D_3+2A_1$ in
$\Delta=B_7$ are reflection groups with Coxeter diagram
\begin{center}
\epsfig{file=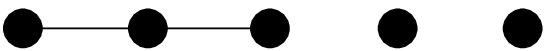,width=0.2\linewidth}
\end{center}
At the same time, there is no automorphism of $G$ taking one of these
groups to another.

It is easy to see that the following is true (cf.~\cite{On}):

\begin{lemma}
\label{b_n-sim}
Given reflection subgroups $H'$ and $H''$ of $G$ of the same type,
there exists an automorphism of $G$ taking $H'$ to $H''$ if and only if
the corresponding root systems $\Delta_1,\,\Delta_2\subset\Delta$ contain
the same number of components of the same type (i.e. we differ components
of the type $A_3$ and $D_3$, $A_1$ and $B_1$, as well as $D_2$, $2A_1$ and
$2B_1$).

\end{lemma}

\begin{col}
\label{b_n-full}
Maximal rank reflection subgroups $H'$ and $H''$ of the same  type
cannot be taken one to another by any  automorphism of $G$ if and only if
the corresponding root systems $\Delta_1,\,\Delta_2\subset\Delta$ are of
the type
\begin{center}
$\Delta_1=B_{l_1}+\dots+B_{l_p}+lB_1+D_{m_1}+\dots+D_{m_r}+mD_2, $
\end{center}
\begin{center}
$
\Delta_2=B_{l_1}+\dots+B_{l_p}+(l+2k)B_1+D_{m_1}+\dots+D_{m_r}+(m-k)D_2,
$
\end{center}
where\quad
$
\sum\limits_{i=1}^p l_i+\sum\limits_{i=1}^r m_i+l+2m=n,
\quad  l_1\ge\dots\ge l_p>1,\quad
m_1\ge\dots\ge m_r>2,\quad k\ne 0.
$

\end{col}

Now we list maximal reflection subgroups of $B_n$.

\begin{lemma}
\label{b_n-max}
Maximal reflection subgroups of the reflection group $B_n$ are either $D_n$
or $B_k+B_{n-k}$, where $k=1,\dots,\left[\frac{n}{2}\right]$.

\end{lemma}

\begin{proof}
Maximality of the groups under consideration follows immediately from
Lemmas~\ref{direct product} and~\ref{b_n}.

Consider a reflection subgroup $H\subset G=B_n$. The
corresponding root system $\Delta_1\subset\Delta$
is of the type described in Lemma~\ref{b_n}, i.e.
$$
\Delta_1=A_{k_1}+\dots+A_{k_s}+B_{l_1}+\dots+B_{l_p}+D_{m_1}+\dots+D_{m_r},
\quad \sum\limits_{i=1}^s(k_i+1)+\sum\limits_{i=1}^p l_i+\sum\limits_{i=1}^r
m_i\le n.
$$

The structure of the root system $B_n$ (see~\cite{Kac}
or~\cite{Bur}) shows that for any component $A_{k_i}$ there exists a root
system of the type $D_{k_i+1}$ which contains $A_{k_i}$ and which is
orthogonal to all the rest components of $\Delta_1$.
Thus,
$$
\Delta_1\subset\Delta_2=B_{l_1}+\dots+B_{l_p}+
D_{k_1+1}+\dots+D_{k_s+1}+D_{m_1}+\dots+D_{m_r},
$$
\begin{center}
$
\sum\limits_{i=1}^s(k_i+1)+\sum\limits_{i=1}^p l_i+\sum\limits_{i=1}^r
m_i\le n.
$
\end{center}

Furthermore, there exists a root system of the type
$D_m$, \ $m={\sum\limits_{i=1}^s(k_i+1)+\sum\limits_{i=1}^r
m_i}$, which contains the system
$D_{k_1+1}+\dots+D_{k_s+1}+D_{m_1}+\dots+D_{m_r}$ and which is orthogonal
to all the components of $\Delta_2$ of the type $B_{l_i}$.
Hence,
$$
\Delta_1\subset\Delta_2\subset\Delta_3=B_{l_1}+\dots+B_{l_p}+D_m,\qquad
\sum\limits_{i=1}^p l_i+m\le n.
$$
Notice also that any component $D_m$ is contained in a root
system $B_m$ which is orthogonal to all components of $\Delta_3$  of
the type $B_{l_i}$. Collecting, in the similar way, all the components of
the type $B_{l_i}$ into root system $B_{l}$ for $l=\sum\limits_{i=1}^p
l_i$, we obtain the claim of the lemma.

\end{proof}

\subsection{Subgroups of $F_4$}
\label{f_4-sec}

Consider the group $G=F_4$ as the Weyl group of the root system
$\Delta=F_4$. Table~10 of~\cite{Dyn} contains the list of all root
subsystems of the root system $F_4$. Now we "forget" lengths of
roots to consider these root subsystems as Coxeter diagrams of
reflection subgroups of the reflection group $F_4$.

\begin{lemma}
\label{f_4}
A group $2B_2$ is a unique maximal rank reflection subgroup of  $F_4$
not appearing in Table~10 of~\cite{Dyn}.

For any two maximal rank reflection subgroups of the same type there
exists some automorphism of $F_4$ taking one of them to another.

\end{lemma}

\begin{proof}

Let  $H\subset G$ be a reflection subgroup of rank 4, such that a Coxeter
diagram of $H$ does not appear in~Table~10 of~\cite{Dyn}, and let
$\Delta_1\subset\Delta$ be the corresponding root system.
Denote by $\alpha_1,\alpha_2,\alpha_3,\alpha_4$ simple roots of
$\Delta_1$.

Since $\Delta_1$ is not a root subsystem of $\Delta$,
there exist short roots $\alpha_i$ and $\alpha_j$, $i=1,\dots,4$,
such that $\alpha_i-\alpha_j\in\Delta$. We may assume that $i=1,j=2$.
In particular, $\alpha_1$ is orthogonal to $\alpha_2$.

Note that a change of lengths of all roots of $\Delta$ induces an
automorphism of the Weyl group $G$. This automorphism takes $H$ to another
maximal rank reflection subgroup $H'$ corresponding to some root system
$\Delta'_1\subset\Delta$; a root system $\Delta'_1$ differs from
$\Delta_1$ by lengths of all roots only. The images
$\alpha'_i\in\Delta'_1$ of the roots $\alpha_i\in\Delta_1$ are simple roots
of $\Delta'_1$.

A Coxeter diagram of $H'$ is obviously the same as one of $H$,
hence $\Delta'_1$ is not a root subsystem of $\Delta$, either.
Since the roots $\alpha'_1$ and $\alpha'_2$ are long, we see that
$\alpha'_3-\alpha'_4\in\Delta$. In particular, $\alpha'_3$ and
$\alpha'_4$ are mutually orthogonal short roots.

Therefore, the set of simple roots of $\Delta'_1$ consists of two short
roots and two long roots, and any two roots of the same length are
mutually orthogonal. Since two roots of different lengths cannot
form an angle $\frac{\pi}{3}$, we obtain that  $H'$ is either of the
type $4A_1$, or $2A_1+B_2$, or $2B_2$, so the first claim of the lemma
is proved.

Further, it is shown in~\cite{Dyn} that any two maximal rank root
subsystems of $\Delta=F_4$ of the same type are conjugated by some
element of the Weyl group $G$. It follows that any two
maximal rank reflection subgroups corresponding to root
systems of the same type are conjugated in $G$.
Extending the inner automorphism group of $G$ by the outer automorphism
described above, we see that for any two maximal rank reflection subgroups
of the same type (different from $2B_2$) there exists an automorphism of
$G$ taking one to another.

Now consider subgroups of the type $2B_2$. Any of them is contained in
some subgroup of the type $B_4$, and for any two subgroups of this type
there exists an automorphism of $G$ taking one of them to another.
Furthermore, any two subgroups of the type $2B_2$ are conjugated
in $B_4$, and any automorphism of $B_4$ may be extended to an automorphism
of $G$.  This proves the second statement of the lemma.

\end{proof}

\begin{lemma}
\label{f_4-max}
A maximal reflection subgroup of $F_4$ is either of the type $B_4$ or
$2A_2$.

\end{lemma}

\begin{proof}
As we have proved before, subgroups of the type $2B_2$ are not maximal.
Thus, any maximal subgroup should either correspond to maximal root subsystem
of $\Delta$ or be a subgroup of smaller rank.  Arguments of the proof
of Lemma~\ref{f_4} applied to a subgroup of smaller rank show
that any such a subgroup is equivalent modulo the automorphism group of
$G$ to a subgroup corresponding to a root subsystem of $\Delta$.  Hence, we
only need to consider subgroups corresponding to maximal root
subsystems.

All the maximal root subsystems of $\Delta=F_4$ are listed
in~Table~12 of~\cite{Dyn} (a subsystem $A_3+A_1$ is not maximal,
see~\cite{On}). Any subgroup of the type $B_3+B_1$ is contained in
some subgroup of the type $B_4$, so it is not maximal, either.
Subgroups of the type $B_4$ are obviously maximal ($[F_4:B_4]=3$).
The same for subgroups of the type $2A_2$ follows from Lemmas~\ref{direct
product} and~\ref{f_4}, and from the description of
reflection subgroups of $B_4$ (Lemma~\ref{b_n}) and $D_4$
(Cor.~\ref{sl} and~\cite[Table 9]{Dyn}).

\end{proof}

\subsection{Subgroups of $G_2$}
\label{g_2-sec}

Rank 2 subgroups of the reflection group $G_2$ are either of the type
$A_2$ or $2A_1$. All the subgroups of the type $A_2$ are conjugated in
$G_2$; any two subgroups of the type $2A_1$ are equivalent modulo the
automorphism group of $G_2$.

\medskip

We present another elementary corollary of results provided
in sections~\ref{b_n-sec}--\ref{g_2-sec} and paper~\cite{Dyn}.

\begin{lemma}
\label{nerazl}
Let $G$ be a finite reflection group different from $H_3$, $H_4$ and
$G_2^{(m)}$ (when $m\ne 2,3,4,6$). Then any two maximal rank
indecomposable reflection subgroups of $G$ of the same type are
equivalent modulo the automorphism group of $G$.

\end{lemma}

\subsection{Subgroups of $G^{(m)}_2$, $H_3$ and $H_4$}
\label{rest-sec}

Rank 2 reflection subgroups of $G_2^{(m)}$ are of the type $G_2^{(d)}$,
for some $d\mid m$. Maximal rank subgroups of the groups $H_3$ and $H_4$
are classified in~\cite{uncut}. In paricular, Theorem~1 of~\cite{uncut}
implies that rank 2 subgroups of $H_3$ are either of the type $G_2^{(5)}$
or $2A_1$. Moreover, it follows from the same theorem  that all 
subgroups of $H_4$ of rank less than 4 are not 
maximal. So, the classification of reflection subgroups of $H_3$ and $H_4$
immediately follows from the classification of reflection subgroups of the
rest finite reflection groups.

\section{Indecomposable subgroups}
\label{simpl}

In this section, we study indecomposable reflection subgroups of
Euclidean reflection groups.

Let $P$ and $F$ be Euclidean Coxeter simplices and $G_P$ be a subgroup
of $G_F$. The group  $G_F$ contains an infinite number of mutually
parallel mirrors. Hence, $P$ may be similar to $F$,
i.e. $\Sigma(P)$ may coincide with $\Sigma(F)$.

\subsection{Similar simplices}

In this section we assume that $\Sigma(F)=\Sigma(P)$.

Let $V$ be a special vertex of $P$.  Without loss of generality we may
assume that $F$ is a fundamental simplex of $G_F$, such that $F$ is
contained in $P$ and contains the vertex $V$. By the
definition of a special vertex, $V$ is a special vertex of $F$.

\begin{lemma}\label{k^n}
1) Let $F$ be a Coxeter simplex in $\E^n$.
For any $k\in \N$ there exists a Coxeter simplex $T$ such that
$\Sigma(T)=\Sigma(F)$,
$G_T\subset G_F$ and $[G_F:G_T]=k^n$.

2) If simplex $F$ is homothetic to $P$ and  $G_P\subset G_F$, then
the dilation factor $k$ is positive integer number and
$[G_F:G_P]=k^n$.

\end{lemma}

\begin{proof}

1) Let $V$ be a special vertex of $F$ and $\mu$ be a facet of $F$
opposite to $V$.
Since $V$ is a special vertex, there exists a mirror $\nu$ containing
$V$ and parallel to $\mu$.

Let $K$ be a simplicial cone with apex $V$ and facets containing
the facets of $F$.
Consider a group $G$ generated by reflections in $\mu$ and $\nu$.
Mirrors of $G$ cut simplices $T_k$ out of $K$. It is clear that
$\Sigma(T_k)=\Sigma(F)$,
$G_{T_k}\subset G_F$ and $[G_F:G_{T_k}]=k^n$.

2) Since $F$ is a fundamental domain of $G_F$, no mirror of $G_F$ parallel
to $\mu$ separates $\mu$ and $\nu$. Thus, simplices $T_k$ are only
simplices contained in $K$ which are homothetic to $F$ and bounded
by mirrors of $G_P$. In other words, $P=T_k$ for some $k$, and
$[G_F:G_P]=k^n$.

\end{proof}

\begin{lemma}
\label{iskl}
Let $F\ne P$ be Coxeter simplices in $\E^n$ such that
$\Sigma(F)=\Sigma(P)$ and $G_P\subset G_F$.
Let $V$ be a common special vertex of $F$ and $P$,
and  $\mu$ be a facet of $P$ opposite to
$V$. The following three conditions are equivalent:

1) $1<[G_F:G_P]<2^n$.

2) No mirror of $G_F$ parallel to $\mu$ intersects the interiour of $P$.

3) One of the following three opportunities holds:

 $\bullet$ $\Sigma(F)=\Sigma(P)=\widetilde C_2$, $[G_F:G_P]=2$;

 $\bullet$ $\Sigma(F)=\Sigma(P)=\widetilde G_2$, $[G_F:G_P]=3$;

 $\bullet$ $\Sigma(F)=\Sigma(P)=\widetilde F_4$, $[G_F:G_P]=4$.

For any of these three cases the subgroup $G_P$ is determined uniquely
up to an automorphism of $G_F$.

\end{lemma}

\begin{proof}

Let $\mu_1,...,\mu_n$ be the facets of $P$ containing $V$.
Consider a diagram
$\Sigma'(F)=\Sigma(F)\setminus v$, where
$v$ is a special node of $\Sigma(F)$
corresponding to $\mu$.
Suppose that
$\Sigma(F)\ne\widetilde  C_2,\widetilde   G_2$ and $\widetilde
F_4$. Then any automorphism of the diagram $\Sigma'(F)$
is a restriction of some automorphism of  the diagram  $\Sigma(F)$.
Thus, the dihedral angles formed up by $\mu$ and
$\mu_1,...,\mu_n$ are uniquely defined.
Hence, $P$ is homothetic to $F$, and condition 2) implies that
$\Sigma(F)=\widetilde   C_2,\widetilde  G_2$ or $\widetilde F_4$.
For any of these cases there exists a unique automorphism of
$\Sigma'(F)$ that is not a restriction of an automorphism of $\Sigma(F)$.
These automorphisms lead to the subgroups shown in Table~\ref{min}.

Therefore, condition 2) implies 3). Now suppose that condition 1)
holds. In this case $P$ is not homothetic to $F$ by the second
part of Lemma~\ref{k^n}.  Thus, $\Sigma(F)=\widetilde   C_2,\widetilde
G_2$ or $\widetilde F_4$, and condition 3) holds.

Evidently, condition 3) implies 1) and 2).

\end{proof}

\begin{table}[!h]
\tcaption{{\bf Three exclusions.}
The table contains all subgroups $G_P\subset G_F$ such that
$\Sigma(F)=\Sigma(P)$ and
$[G_F:G_P]<2^n$.
Vectors $\xi_1,...,\xi_{n+1}$ are the outward normals to the facets
$f_i$  of $F$. The facets are indexed as it is shown in the second
column.\\
\phantom{ww}By $r_i$ we denote the reflection with respect to $f_i$.
The normals to the facets of  $P$ are expressed in terms of 
$\xi_1,...,\xi_{n+1}$\vphantom{$\int\limits_a$}.
}
\label{min}
\begin{center}
\begin{tabular}{|c|@{\quad}c@{\quad}|@{\quad}c@{\quad}|c|c|}
\hline
Type&$\Sigma(F)$&$\Sigma(P)$&
\begin{tabular}{c}
additional\\
\raisebox{4pt}[1pt][0pt]{vector}
\end{tabular}
&index\\
\hline
$\widetilde C_2$
&\raisebox{-7pt}[10pt][12pt]{
\psfrag{1}{\scriptsize $1$}
\psfrag{2}{\scriptsize $2$}
\psfrag{3}{\scriptsize $3$}
\epsfig{file=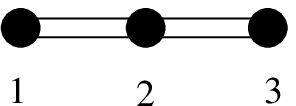,width=0.1\linewidth}
}
&\raisebox{-7pt}[10pt][12pt]{
\psfrag{1}{\scriptsize $2$}
\psfrag{2}{\scriptsize $3$}
\psfrag{3}{\scriptsize $4$}
\epsfig{file=pic/tc_2.eps,width=0.1\linewidth}
}
&$\xi_4=r_1(\xi_2)$
&2\\
\hline
$\widetilde G_2$
&\raisebox{-7pt}[10pt][12pt]{
\psfrag{1}{\scriptsize $1$}
\psfrag{2}{\scriptsize $2$}
\psfrag{3}{\scriptsize $3$}
\epsfig{file=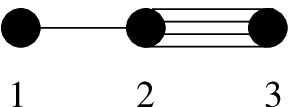,width=0.1\linewidth}
}
&\raisebox{-7pt}[10pt][12pt]{
\psfrag{1}{\scriptsize $4$}
\psfrag{2}{\scriptsize $3$}
\psfrag{3}{\scriptsize $2$}
\epsfig{file=pic/tg_2.eps,width=0.1\linewidth}
}
&
$\xi_4=r_1r_2(\xi_3)$
&3\\
\hline
$\widetilde F_4$
&\raisebox{-7pt}[10pt][12pt]{
\psfrag{1}{\scriptsize $1$}
\psfrag{2}{\scriptsize $2$}
\psfrag{3}{\scriptsize $3$}
\psfrag{4}{\scriptsize $4$}
\psfrag{5}{\scriptsize $5$}
\epsfig{file=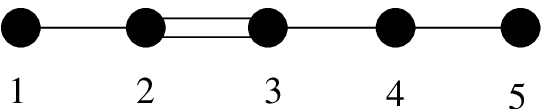,width=0.18\linewidth}
}
&\raisebox{-7pt}[10pt][12pt]{
\psfrag{1}{\scriptsize $$4}
\psfrag{2}{\scriptsize $3$}
\psfrag{3}{\scriptsize $2$}
\psfrag{4}{\scriptsize $1$}
\psfrag{5}{\scriptsize $6$}
\epsfig{file=pic/tf_4.eps,width=0.18\linewidth}
}
&$\xi_6=r_5r_4r_3(\xi_2)$
&4\\
\hline
\end{tabular}
\end{center}

\end{table}


\begin{lemma}
\label{uniqueP-P}
Let $P$ and $F$ be similar Euclidean Coxeter simplices such that
$G_P\subset G_F$. Then

1) The index $[G_F:G_P]$ determines the subgroup
$G_P\subset G_F$ by a unique way up to an automorphism of $G_F$.

2) If the subgroup $G_P\subset G_F$ is maximal then either $F$ is
homothetic to $P$, or $G_P\subset G_F$ is one of the subgroups presented in
Table~\ref{min}.

\end{lemma}

\begin{proof}
Let $V$ be a common special vertex of $P$
and $F$, and assume that $P$ contains $F$.

Suppose that $\Sigma(F)\ne\widetilde  C_2,\widetilde  G_2$ and
$\widetilde F_4$. Then we can assume that $P=h(F)$, where $h$
is a homothety centered at $V$ (see Lemma~\ref{iskl}), and both
statements of the lemma are evident.

Now suppose that $\Sigma(F)=\widetilde  C_2,\ \widetilde  G_2$ or
$\widetilde F_4$.
 Let $\mu_P$ be the facet of $P$ opposite to
 $V$.  Let $\mu$ be the closest to $V$ mirror of $G_F$ parallel to $\mu_P$
 such that $V\notin \mu$ and $\mu$ intersects $P$ ($\mu$ may coincide with
$\mu_P$).
The mirror $\mu$ cuts some simplex $T$ out of $P$, $T=h(P)$,
 where $h$ is a homothety centered at $V$ with dilation factor $\frac{1}{k}$,
 $k\in \N$ ($k$ may equal one).  Consider the subgroup
 $G_T\subset G_F$.  Notice that $G_P\subset G_T$, since $V$ is a special
vertex, and $\mu$ is the closest mirror to $V$ parallel to $\mu_P$ and not
containing $V$. By Lemma~\ref{iskl}, either $T=F$ or $[G_F:G_T]=2,3,4$
respectively for $\Sigma(F)=\widetilde C_2,\widetilde G_2$ and $\widetilde
F_4$.  Thus, $[G_F:G_P]=mk^n$, where $k\in \N$, and either $m=1$ or
$m=2,3$ and $4$ respectively for $\Sigma(F)=\widetilde C_2,\widetilde G_2$
and $\widetilde F_4$.  Clearly, such a subgroup $G_P\subset G_F$ is
completely determined by $m$ and $k$, and the numbers $m$ and $k$ are
determined uniquely by the index $[G_F:G_P]=mk^n$. This proves the
 first part of the lemma. If the subgroup $G_P\subset G_F$ is maximal then
 either $m$ or $k$ is equal to one. This implies the second statement.

\end{proof}

We provide another one fact concerning similar simplices in $\E^n$.

\begin{lemma}
\label{easy_to_check}
Let $F$ and ${P}$ be simplices in $\E^n$, where $\Sigma(F)=\Sigma(P)$
and $G_{P}\subset G_F$.
Then any automorphism of $G_{P}$ is a restriction of some
automorphism of $G_{F}$.

\end{lemma}

\begin{proof}

 Let $\varphi$ be an automorphism of  $G_{P}$. It is well-known (see,
 e.g.,~\cite{ChD}) that $\varphi$ takes reflections to reflections. This
 means that we may consider $\varphi$ as an isometry of $\E^n$ preserving
 the set of mirrors of $G_P$.

  Let $V$ be a special vertex of $P$, and $F_1$ be
 a fundamental simplex of $G_F$ that contains the vertex $V$ and
 lies inside of $P$.  Such a simplex $F_1$ is unique since
 $\Sigma(F)=\Sigma(P)$. Then $\varphi(V)$ is a special vertex of
 a fundamental simplex $\varphi(P)$ of $G_P$, and $\varphi(F_1)$ is a
 simplex congruent to $F$ contained in $\varphi(P)$ and containing the
 vertex $\varphi(V)$.  Denote by $F_2$ the fundamental simplex of $G_F$
 contained in $\varphi(P)$ and containing $\varphi(V)$.

 If
 $\Sigma(F)\ne \widetilde C_2, \widetilde G_2$ and $\widetilde F_4$,
 then  there exists a unique simplex of given volume with Coxeter
 diagram $\Sigma(F)$ containing $\varphi(V)$ as a special vertex
 and contained in $\varphi(P)$. Hence, $F_2=\varphi(F_1)$,
 so $\varphi$  is an automorphism of $G_F$.

 Suppose now that
 $\Sigma(F)= \widetilde C_2, \widetilde G_2$ or $\widetilde
 F_4$. If $\varphi(F_1)=F_2$ we have nothing to prove.
 Suppose  that $\varphi(F_1)\ne F_2$.
 Then either $F_1$ is homothetic to $P$, $F_2$ is not homothetic to
 $\varphi(P)$,  or  $F_1$ is not homothetic to $P$, $F_2$ is  homothetic to
 $\varphi(P)$. Thus, one of the indices
 $[G_{F_1}:G_{P}]$ and $[G_{F_2}:G_{\varphi(P)}]$ equals $k^n$,
 and another one equals $2k^n$,  $3k^n$ and $4k^n$
 respectively for the cases $\Sigma(F)= \widetilde C_2, \widetilde
  G_2$ and $\widetilde  F_4$. This contradicts to the fact that
 $G_{F_1}=G_{F_2}$ and $G_{P}=G_{\varphi(P)}$.

\end{proof}

\subsection{Simplices with distinct Coxeter diagrams}

Now suppose that $P$ is not similar to $F$.

Let $\mu(v)$ denote the facet of Coxeter simplex $T$ corresponding to a
 node $v$ of $\Sigma(T)$, and let $V$ be the vertex of $T$ opposite to
 $\mu(v)$.  Denote by $G_{T\setminus \mu(v)}$ the subgroup of $G_T$
generated by reflections in the facets of $T$
different from  $\mu(v)$.
In other words,  $G_{T\setminus \mu(v)}$ is the stabilizer of $V$ in $G_T$.

\begin{lemma}
\label{sub}
Let $F$ and $P$ be Coxeter simplices in $\E^n$, where $G_P\subset G_F$.
Then for any node $v$ of $\Sigma(P)$ there exists a node
$w$ of  $\Sigma(F)$ such that $G_{P\setminus \mu(v)}
\subset G_{F_w  \setminus \mu(w)}$ for some fundamental simplex $F_w$ of
$G_F$.

\end{lemma}

\begin{proof}
Let $V$ be the vertex of $P$ opposite to $\mu(v)$. Let $F_w$ be a
fundamental simplex of $G_F$ such that $F_w$ is contained in $P$ and
contains $V$. Let $w$ be the node of the diagram $\Sigma(F)$ corresponding
to the facet of $F_w$ opposite to $V$. The stabilizer of $V$ in $G_P$ is a
subgroup of the  stabilizer of $V$ in $G_F$. At the same time, the
stabilizers of $V$ in $G_P$ and $G_F$ coincide with
$G_{P\setminus \mu(v)}$ and $G_{F_w  \setminus \mu(w)}$
respectively, so the lemma is proved.

\end{proof}

Lemma~\ref{sub} gives a necessary condition for $G_P$ to be a subgroup of
$G_F$ in terms of Coxeter diagrams $\Sigma(F)$ and $\Sigma(P)$. This
condition is easy to check:  the groups $G_{P\setminus \mu(v)}$ and $
G_{g(F) \setminus \mu(u)}$ act in the spherical space $\Sp^{n-1}$, and
we can use results of Section~\ref{sph}. Straightforward check of the
condition described in Lemma~\ref{sub} gives rise to the following claim:

\begin{lemma}
\label{list}
Let $F$ and $P$ be Coxeter simplices in $\E^n$ such that $G_P\subset G_F$
and $\Sigma(P)\ne\Sigma(F)$. Then the pair of diagrams
$(\Sigma(P),\Sigma(F))$ coincides with one of the pairs listed in
Table~\ref{res}.

\end{lemma}

In particular, Lemma~\ref{sub} implies that
if $G_{P}\subset G_{F}$ and
$\Sigma(P)= \widetilde  C_2,\widetilde  G_2$ or $\widetilde F_4$,
then $\Sigma(F)=\Sigma(P)$.

Further, for any pair $(F,P)$ with $\Sigma(P)\ne\Sigma(F)$
satisfying Lemma~\ref{sub} we find simplices $F_1$ and $P_1$ such that
$\Sigma(F_1)=\Sigma(F)$, $\Sigma(P_1)=\Sigma(P)$ and $G_{P_1}\subset
G_{F_1}$. Examples of such simplices are presented in the right column of
Table~\ref{res}.

\begin{lemma}
\label{spec}
Let $F$ and $P$ be Coxeter simplices in $\E^n$ such that the subgroup
$G_P\subset G_F$ is maximal and $\Sigma(P)\ne
\Sigma(F)$.
Then there exist a special vertex  $V$ of $P$ and
a fundamental simplex  $F_1$ of $G_F$ such that
$V$ is a special vertex of $F_1$.

\end{lemma}

\begin{proof}

Let $L$ be the set of all mirrors of
$G_F$ parallel to mirrors of $G_P$.  Denote by $G_L$
the group generated by reflections in all mirrors
contained in $L$.  Obviously, $G_{P}\subseteq G_{L}\subseteq G_F$ and the
subgroup $G_{P}\subset G_F$ is maximal. Thus, either $G_L=G_P$ or
$G_L=G_F$.

Let $F_1$ be a fundamental simplex of $G_F$ contained in $P$, and let
$V$ be a special vertex of $F_1$.
For any mirror $\mu$ of $G_P$ there exists a mirror
containing $V$ and parallel to $\mu$.
If $G_L=G_P$ then any of these mirrors is contained in $G_P$, so $V$
is a special  vertex of $P$.

Now, suppose that $G_L=G_F$. Then the set of linear parts of elements of
$G_P$ coincides with the set of linear parts of elements of $G_F$.  It
follows that either $\Sigma(F)=\widetilde C_n$, $\Sigma(P)=\widetilde B_n$
or $\Sigma(F)=\widetilde B_n$, $\Sigma(P)=\widetilde C_n$.  In both cases
each special vertex of $P$ is a special vertex of some fundamental simplex
of $G_F$.

\end{proof}

\begin{lemma}
\label{maxi}

Let $G_P\subset G_F$ be a maximal subgroup, where $F$ and $P$ are
Euclidean  Coxeter simplices with $\Sigma(P)\ne \Sigma(F)$. Then $G_P$
is determined by the pair $(\Sigma(P), \Sigma(F))$
uniquely up to an automorphism of $G_F$.
\end{lemma}

\begin{proof}
Let $P_1$ and $P_2$ be simplices similar to  $P$,
and suppose the subgroups $G_{P_i}\subset G_F$ to be maximal.
It is sufficient to show that there exists an automorphism of $G_F$
taking $P_1$ to $P_2$.

By Lemma~\ref{spec}, there exist a special vertex $V_1$ of $P_1$
and a fundamental simplex $F_1$ of $G_F$ such that $V_1$ is a special
vertex of $F_1$.
Let $V_2$ be a special vertex of $P_2$, and $F_2$ be
a fundamental simplex of $G_F$ such that  $V_2$ is a special
vertex of $F_2$.

Consider an automorphism $\varphi$ of $G_F$ taking
$F_1$ to $F_2$.
Notice that for any Euclidean Coxeter simplex $T$ the following holds:
for any two special vertices $V$ and $U$ of $T$ there exists a symmetry
of $T$ exchanging $U$ and $V$.
Hence, there exists an automorphism $\psi$ of $G_F$ taking $F_2$ to
itself and satisfying $\psi\circ\varphi(V_1)=V_2$.
Denote $\psi\circ\varphi(P_1)$ by $P_3$.

Let $K_2$ and $K_3$ be the minimal cones with an apex $V_2$ containing
$P_2$ and $P_3$ respectively.

Consider the stabilizers $G'_F(V_2)$ and  $G'_P(V_2)$ of $V_2$
in $G_F$ and $G_P$ respectively.
As it is shown in Section~\ref{sph} (Lemma~\ref{nerazl}), an
indecomposable maximal rank finite subgroup $G'_{P}(V_2)\subset
G'_{F}(V_2)$ is determined $G'_{F}(V_2)$ by $\Sigma'(P)$  and
$\Sigma'(F)$ uniquely up to an automorphism of 
$G'_{F}(V_2)$.  Hence, there exists an
automorphism $\rho'$ of $G'_F(V_2)$ sending $K_3$ to $K_2$. Since
$\Sigma(P)\ne\Sigma(F)$, we may assume that
$\Sigma(P)\ne  \widetilde  C_2,\widetilde  G_2$ and $\widetilde F_4$.
Therefore, $\rho'$ is a restriction of some automorphism $\rho$ of $G_F$.
Let $P_4=\rho(P_3)$.

Since  $\Sigma(P)\ne  \widetilde  C_2,\widetilde  G_2$ and $\widetilde
F_4$, any two mirrors cutting a simplex similar to $P$ out of $K_2$ are
mutually parallel. The subgroup $G_{P_1}\subset G_F$ is maximal,
consequently $P_2$ is cut off  by a closest to $V$ mirror
described above. The same is true for $P_4$.  Thus,
$P_4=P_2$, so $\rho\circ\psi\circ\varphi(P_1)=P_2$.

\end{proof}

\begin{lemma}
\label{hom}
Let $G_{P}\subset G_{P_1}\subset G_F$, where  $F$  and $P$ are
Euclidean Coxeter simplices with a common vertex $V$. Suppose that
$F$ is contained in $P$, and
let $P$ be the image of $P_1$ under the homothety centered at $V$
with dilation factor $k\in \N$.
Then there exists a homothety with
dilation factor $k$ taking $F$ to $F_1$,  such that $G_{P}\subset
G_{F_1}\subset G_F$.

\end{lemma}

The proof is evident: take $V$ as the center of homothety.

\vspace{8pt}

Summing up the above, we obtain the following

\begin{theorem}
Let $F$ and $P$ be Coxeter simplices in $\E^n$, and $G_P\subset G_F$.
Then there exists a sequence of subgroups
$G_P=G_{F_l}\subset G_{F_{l-1}}\subset\dots\subset G_{F_1}\subset
G_{F_0}=G_F$, where $G_{F_{i+1}}\subset G_{F_i}$
is a subgroup described either in Table~\ref{min},
or in Table~\ref{res}, or in Lemma~\ref{k^n}.

The subgroup $G_P\subset G_F$ is determined by the index
$[G_F:G_P]$ uniquely up to an automorphism of $G_F$.

\end{theorem}

\begin{proof}

Since $[G_F: G_P]<\infty$, there exists  a sequence  $G_P=G_{F_l}\subset
G_{F_{l-1}}\subset\dots\subset G_{F_1}\subset G_{F_0}=G_F$ such that
any subgroup $G_{F_{i+1}}\subset G_{F_i}$ is maximal. Consider those parts
of the sequence, for which $\Sigma(F_{i+1})=\Sigma(F_i)$.
By Lemma~\ref{uniqueP-P}, a subgroup  $G_{F_{i+1}}\subset
G_{F_i}$ is one described either in Table~\ref{min}, or in
Lemma~\ref{k^n}. Now consider those parts for which
$\Sigma(F_{i+1})\ne\Sigma(F_i)$.  By Lemmas~\ref{list}
and~\ref{maxi}, these parts are described in Table~\ref{res},
and the first statement of the theorem is proved.

We only left to show that  the index $[G_F:G_P]$ determines
the subgroup $G_P\subset G_F$ uniquely up to an automorphism of
$G_F$.

 We say that $\Sigma_q\subset\Sigma_{q-1}\subset \dots
\subset\Sigma_1$
is an {\it admissible sequence} of diagrams for the subgroup
$G_P\subset G_F$ if there exist simplices
$T_q=P$, $T_{q-1}$, $\dots$, $T_1$, $T_0=F$ satisfying the following
conditions:

(1) $\Sigma(T_i)=\Sigma_i, 1\le i\le q$,

(2) $\Sigma_1=\Sigma(F)$,

(3) $G_{T_{i+1}}\subset G_{T_{i}}$ is a maximal subgroup if
$\Sigma(T_i)\ne\Sigma(T_{i+1})$,

(4) $\Sigma(T_i)=\Sigma(T_j),\ i<j,$  if and only if
$i=0,\,j=1$.\\
The sequence of subgroups
$G_P=G_{T_q}\subset G_{T_{q-1}}\subset\dots\subset  G_{T_1}\subset
G_{T_0}=G_F$ satisfying conditions
(1)--(4) we also call {\it admissible}.

Now we will show that for any subgroup $G_P\subset G_F$ there exists
an admissible sequence of subgroups.

Let
$G_P=G_{F_l}\subset G_{F_{l-1}}\subset\dots\subset G_{F_1}\subset
G_{F_0}=G_F$ be a sequence of subgroups described in the theorem.
Suppose that $\Sigma(F_j)= \widetilde  C_2,\widetilde  G_2$, or
$\widetilde F_4$, and $j>0$. Then
$\Sigma(F)=\Sigma(F_1)=\dots=\Sigma(F_j)$, so the subgroup
$G_{F_j}\subset  G_{F_0}$ is determined by the index
$[G_{F_0}:G_{F_j}]$ by Lemma~\ref{uniqueP-P}.
Hence, we may assume that
$\Sigma(F_j)\ne \widetilde  C_2,\widetilde  G_2,\widetilde F_4$ if
$j>1$.

Further, let $\Sigma(F_i)=\Sigma(F_j),\ i<j,\ 1<j$.
Using Lemma~\ref{hom}, we can subtract the subgroups
$G_{F_{i+1}},\dots,G_{F_{j}}$ from the sequence in the following way:
if $\Sigma(F)=\Sigma(F_1)$ then apply a homothety with dilation factor
 $[G_{F_{i}}:G_{F_j}]$  to the simplices $F_1,\dots,F_{i}$;
if $\Sigma(F)\ne \Sigma(F_1)$ then apply a homothety with factor
 $[G_{F_{i}}:G_{F_j}]$ to the simplices $F_0,\dots,F_{i}$
and insert obtained subgroups between  $G_{F}$ and $G_{F_{j+1}}$.
Note that after any of these procedures simplex $T_1$ is similar to
$F$.

Thus, we need at most $l-1$ steps to transform the initial sequence to the
required one (if no simplex in the initial sequence is similar to $F$
we simply insert $T_1=F$ between $F$ and $F_1$).
Now, using Table~\ref{res} and Lemma~\ref{maxi},
it is easy to find all admissible sequences of diagrams
$\Sigma_q\subset\Sigma_{q-1}\subset \dots
\subset\Sigma_1$.
The sequences with  $q\ge 3$ are listed below
(for the case $q=2$ see Table~\ref{res}).

\vspace{8pt}
\begin{center}
\begin{tabular}{|c|c|}
\hline
$\Sigma_q\subset\dots\subset\Sigma_1$  & $[G_{T_1}:G_{P}]$
\vphantom{$\int\limits_A^A$}\\
\hline
$\widetilde D_4\subset \widetilde B_4 \subset \widetilde F_4$ &
$2\cdot 3$\vphantom{$\int\limits^A$}\\
$\widetilde C_4\subset \widetilde B_4 \subset \widetilde F_4$ &
$2^3\cdot 3$\vphantom{$\int\limits^a$}\\
$\widetilde D_n\subset \widetilde B_n \subset \widetilde C_n$ &
4\vphantom{$\int\limits_A^a$}\\
\hline
\end{tabular}
\end{center}
\vspace{8pt}

For each subgroup we found some admissible sequence of
diagrams.
Notice that for any pair $(\Sigma_q,\Sigma_1$) there exists at most one
admissible sequence $\Sigma_q\subset\Sigma_{q-1}\subset \dots
\subset\Sigma_1$.
Hence, it is sufficient to show that each admissible sequence of
diagrams corresponds to at most one subgroup of given index.

We are rest with two cases:  $q=2$ and $q=3$
(see Lemma~\ref{uniqueP-P} for the case $q=1$).

Suppose that $q=2$.
By Lemma~\ref{maxi}, the subgroup
$G_P\subset G_{T_1}$ is  determined uniquely up to an automorphism of
$G_{T_1}$.
In particular, $[G_{T_1}:G_{P}]$ is uniquely determined.
Hence, the index $[G_F:G_{T_1}]$ is determined, too.
By Lemma~\ref{uniqueP-P}, the subgroup $G_{T_1}\subset G_F$
is determined up to an automorphism of $G_F$.
Since $\Sigma(F)=\Sigma(T_1)$,  Lemma~\ref{easy_to_check}
implies that any automorphism of
$G_{T_1}$ is a restriction of some automorphism of $G_F$.
Thus, $G_P\subset G_F$ is uniquely determined by
$[G_{F}:G_{P}]$.

Suppose that $q=3$. A direct examination
shows that for each of three admissible
sequences $\Sigma_3\subset\Sigma_2\subset\Sigma_1$ there exists a
unique (up to an automorphism of $G_{T_1}$)
subgroup $G_{T_3}\subset G_{T_1}$.
Applying the arguments used for the case of $q=2$, we obtain the theorem.

\end{proof}


\begin{table}[!h]
\tcaption{{\bf Indecomposable maximal subgroups.} 
Finite index indecomposable maximal reflection
subgroups of Euclidean reflection groups are listed in the table.
Notation is the same as in Table~\ref{min}.}
\label{res}
\bigskip
\begin{center}
\begin{tabular}{|@{}c@{}c@{\quad }|@{}l@{}|c|}
\hline
\multicolumn{2}{|c|}{$\Sigma(F)$}&\multicolumn{1}{c|}{$\Sigma(P)$}
&index\\
\hline
\begin{tabular}{c}
{$ {\widetilde B_n}$}\\
$n\ge 3$
\end{tabular}
&\raisebox{-13pt}[5pt][5pt]{
\psfrag{1}{\scriptsize $1$}
\psfrag{2}{\scriptsize $3$}
\psfrag{3}{\scriptsize $4$}
\psfrag{4}{\scriptsize ${n}$}
\psfrag{5}{\scriptsize ${n+1}$}
\psfrag{6}{\scriptsize $2$}
\epsfig{file=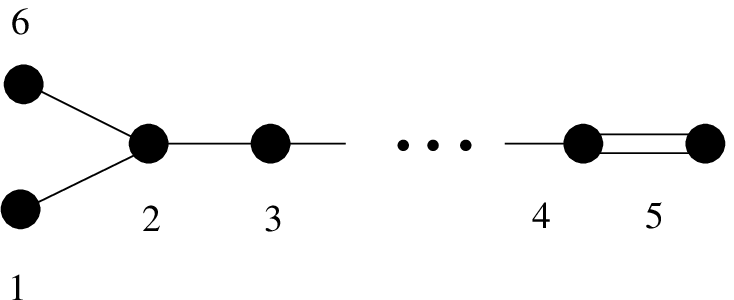,width=0.2\linewidth}
}
&\begin{tabular}{lc}
{$ {\widetilde D_n}$}
&\raisebox{-8pt}[35pt][7pt]{
\psfrag{1}{\scriptsize $1$}
\psfrag{2}{\scriptsize $3$}
\psfrag{3}{\scriptsize $4$}
\psfrag{4}{\scriptsize $\!\!{n-1}$}
\psfrag{5}{\scriptsize $n$}
\psfrag{6}{\scriptsize ${n+2}$}
\psfrag{7}{\scriptsize $2$}
\epsfig{file=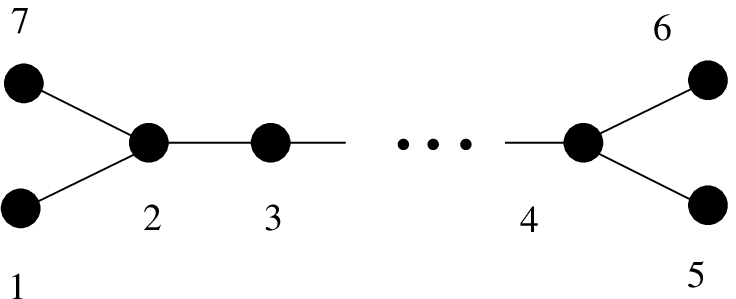,width=0.2\linewidth}
}
\\
\multicolumn{2}{c}{\phantom{$rr_3r_3r_3r_3r_3r_3r_3r_3r_3r_3$}$\xi_{n+2}=r_{n+1}(\xi_n)$\phantom{$r_3r_3r_3r_3r_3r_3r_3r$}\vphantom{$\int\limits_a$}}
\end{tabular}
&2\\
\hline
\begin{tabular}{c}
{$ {\widetilde C_n}$}\\
$n\ge 3$
\end{tabular}
&\raisebox{-4pt}[10pt][5pt]{
\psfrag{1}{\scriptsize $1$}
\psfrag{2}{\scriptsize $2$}
\psfrag{3}{\scriptsize $3$}
\psfrag{4}{\scriptsize ${n}$}
\psfrag{5}{\scriptsize ${n+1}$}
\epsfig{file=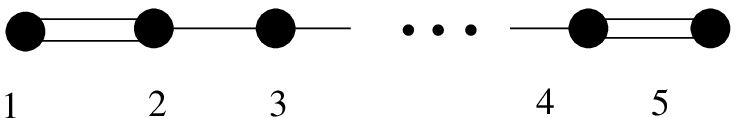,width=0.2\linewidth}
}
&\begin{tabular}{cc}
{$ {\widetilde B_n}$}
&\raisebox{-8pt}[35pt][7pt]{
\psfrag{1}{\scriptsize $2$}
\psfrag{2}{\scriptsize $3$}
\psfrag{3}{\scriptsize $4$}
\psfrag{4}{\scriptsize ${n}$}
\psfrag{5}{\scriptsize ${n+1}$}
\psfrag{6}{\scriptsize ${n+2}$}
\phantom{$r_3r$}%
\epsfig{file=pic/tb_n.eps,width=0.2\linewidth}
}
\\
\multicolumn{2}{c}{\phantom{$r_3r_3r_3r_3r_3r_3r_3r_3r_3r_3r_3$}$\xi_{n+2}=r_{1}(\xi_2)$\phantom{$r_3r_3r_3r_3r_3r_3r_3$}\vphantom{$\int\limits_a$}}
\end{tabular}
&2\\
\hline
\begin{tabular}{c}
{$ {\widetilde B_n}$}\\
$n\ge 3$
\end{tabular}
&\raisebox{-12pt}[30pt][7pt]{
\psfrag{1}{\scriptsize $1$}
\psfrag{2}{\scriptsize $3$}
\psfrag{3}{\scriptsize $4$}
\psfrag{4}{\scriptsize ${n}$}
\psfrag{5}{\scriptsize ${n+1}$}
\psfrag{6}{\scriptsize $2$}
\epsfig{file=pic/tb_n.eps,width=0.2\linewidth}
}
&\begin{tabular}{cc}
{$ {\widetilde C_n}$}
&
\psfrag{1}{\scriptsize $\!\!\!\!\!{n+2}$}
\psfrag{2}{\scriptsize $1$}
\psfrag{3}{\scriptsize $3$}
\psfrag{4}{\scriptsize ${n}$}
\psfrag{5}{\scriptsize ${n+1}$}
\phantom{$r_3$}%
\epsfig{file=pic/tc_n.eps,width=0.2\linewidth}
\\
\multicolumn{2}{c}{\phantom{$r_3r_3r_3r_3r_3r_3r_3r_3r_3$}$\xi_{n+2}=r_{2}r_3r_4\dots r_n(\xi_{n+1})$\phantom{$r_3r_3r_3r_3$}\vphantom{$\int\limits_a$}}
\end{tabular}
&$2^{n-1}$\\
\hline
{$ {\widetilde E_8}$}
&\raisebox{-11pt}{
\psfrag{1}{\scriptsize $1$}
\psfrag{2}{\scriptsize $2$}
\psfrag{3}{\scriptsize $3$}
\psfrag{4}{\scriptsize ${4}$}
\psfrag{5}{\scriptsize ${5}$}
\psfrag{6}{\scriptsize $6$}
\psfrag{7}{\scriptsize ${7}$}
\psfrag{8}{\scriptsize ${8}$}
\psfrag{9}{\scriptsize $9$}
\epsfig{file=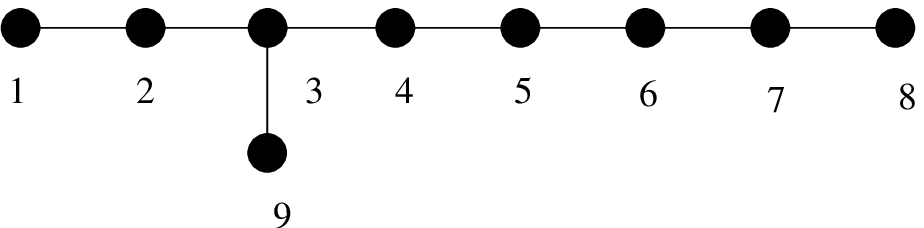,width=0.2\linewidth}
}
&\begin{tabular}{cc}
{$ {\widetilde D_8}$}
&\raisebox{-5pt}[30pt][5pt]{
\psfrag{2}{\scriptsize $2$}
\psfrag{3}{\scriptsize $3$}
\psfrag{4}{\scriptsize ${4}$}
\psfrag{5}{\scriptsize ${5}$}
\psfrag{6}{\scriptsize $6$}
\psfrag{7}{\scriptsize ${7}$}
\psfrag{8}{\scriptsize ${8}$}
\psfrag{9}{\scriptsize $9$}
\psfrag{10}{\scriptsize ${10}$}
\phantom{$r_3$}%
\epsfig{file=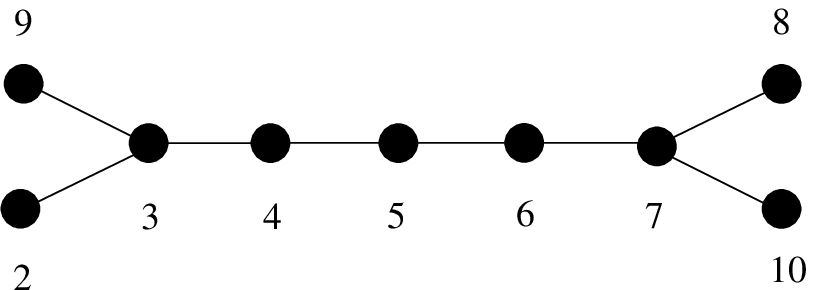,width=0.2\linewidth}
}
\\
\multicolumn{2}{c}{\phantom{$r_3r_3r_3$}$\xi_{10}=r_1r_2r_3r_4r_5r_6r_9r_3r_4r_5r_2r_1r_3r_9r_2r_4(\xi_3)$\phantom{$r_3r_3$}\vphantom{$\int\limits_a$}}
\end{tabular}
&$2\cdot 3^3\cdot 5$\\
\hline
{$ {\widetilde E_8}$}
&\raisebox{-11pt}{
\psfrag{1}{\scriptsize $1$}
\psfrag{2}{\scriptsize $2$}
\psfrag{3}{\scriptsize $3$}
\psfrag{4}{\scriptsize ${4}$}
\psfrag{5}{\scriptsize ${5}$}
\psfrag{6}{\scriptsize $6$}
\psfrag{7}{\scriptsize ${7}$}
\psfrag{8}{\scriptsize ${8}$}
\psfrag{9}{\scriptsize $9$}
\epsfig{file=pic/te_8.eps,width=0.2\linewidth}
}
&\begin{tabular}{cc}
{$ {\widetilde A_8}$}
&\raisebox{-7pt}[40pt][5pt]{
\psfrag{1}{\scriptsize $1$}
\psfrag{2}{\scriptsize $2$}
\psfrag{3}{\scriptsize $7$}
\psfrag{4}{\scriptsize ${8}$}
\psfrag{5}{\scriptsize ${10}$}
\phantom{$r_3$}%
\epsfig{file=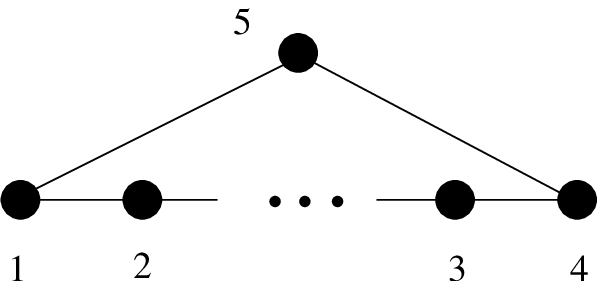,width=0.2\linewidth}
}
\\
\multicolumn{2}{c}{$\xi_{10}=r_9r_3r_4r_5r_6r_7r_2r_3r_4r_5r_6r_9r_3r_4r_5r_2r_1r_3r_9r_2r_4(\xi_3)$\vphantom{$\int\limits_a$}} 
\end{tabular}
&$2^7\cdot 3^2\cdot 5$\\
\hline
{$ {\widetilde E_7}$}
&\raisebox{-13pt}{
\psfrag{1}{\scriptsize $1$}
\psfrag{2}{\scriptsize $2$}
\psfrag{3}{\scriptsize $3$}
\psfrag{4}{\scriptsize ${4}$}
\psfrag{5}{\scriptsize ${5}$}
\psfrag{6}{\scriptsize $6$}
\psfrag{7}{\scriptsize ${7}$}
\psfrag{8}{\scriptsize ${8}$}
\epsfig{file=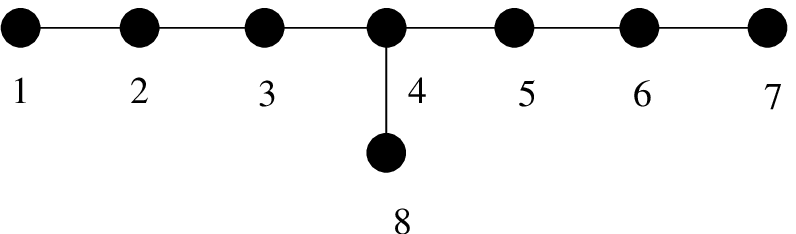,width=0.2\linewidth}
}
&\begin{tabular}{cc}
{$ {\widetilde A_7}$}
&\raisebox{-7pt}[38pt][5pt]{
\psfrag{1}{\scriptsize $1$}
\psfrag{2}{\scriptsize $2$}
\psfrag{3}{\scriptsize $6$}
\psfrag{4}{\scriptsize ${7}$}
\psfrag{5}{\scriptsize ${9}$}
\phantom{$r_3r_3$}%
\epsfig{file=pic/ta_n.eps,width=0.2\linewidth}
}
\\
\multicolumn{2}{c}{\phantom{$r_3r_3r_3r_3r_3r_3r_3r_3$}$\xi_{9}=r_8r_4r_3r_2r_5r_6r_4r_8r_3r_5(\xi_4)$\phantom{$r_3r_3r_3$}\vphantom{$\int\limits_a$}}
\end{tabular}
&$2^4\cdot 3^2$\\
\hline
{$ {\widetilde F_4}$}
&\raisebox{-2pt}{
\psfrag{1}{\scriptsize $1$}
\psfrag{2}{\scriptsize $2$}
\psfrag{3}{\scriptsize $3$}
\psfrag{4}{\scriptsize ${4}$}
\psfrag{5}{\scriptsize ${5}$}
\epsfig{file=pic/tf_4.eps,width=0.14\linewidth}
}
&\begin{tabular}{cc}
{$ {\widetilde B_4}$}
&\raisebox{-2pt}[20pt][5pt]{
\psfrag{1}{\scriptsize $2$}
\psfrag{2}{\scriptsize $3$}
\psfrag{3}{\scriptsize ${4}$}
\psfrag{4}{\scriptsize $5$}
\psfrag{5}{\scriptsize ${6}$}
\phantom{$r_3r_3$}%
\epsfig{file=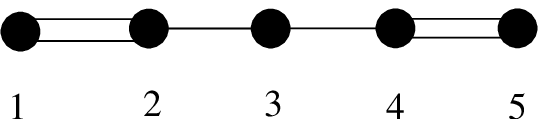,width=0.14\linewidth}
}
\\
\multicolumn{2}{c}{\phantom{$r_3r_3r_3r_3r_3r_3r_3r_3r_3r_3r_3r_3$}$\xi_{6}=r_1r_2(\xi_3)$\phantom{$r_3r_3r_3r_3r_3r_3$}\vphantom{$\int\limits_a$}}
\end{tabular}
&3\\
\hline
{$ {\widetilde G_2}$}
&\raisebox{-2pt}{
\psfrag{1}{\scriptsize $1$}
\psfrag{2}{\scriptsize $2$}
\psfrag{3}{\scriptsize $3$}
\epsfig{file=pic/tg_2.eps,width=0.07\linewidth}
}
&\begin{tabular}{cc}
\raisebox{-7pt}{$ {\widetilde A_2}$}
&\raisebox{-17pt}[33pt][17pt]{
\psfrag{1}{\scriptsize $4$}
\psfrag{2}{\scriptsize $3$}
\psfrag{3}{\scriptsize $2$}
\phantom{$r_3r_3$}%
\epsfig{file=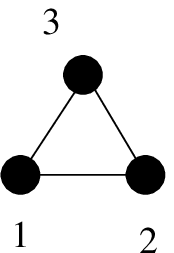,width=0.07\linewidth}
}
\\
\multicolumn{2}{c}{\phantom{$r_3r_3r_3r_3r_3r_3r_3r_3r_3r_3r_3r_3$}$\xi_4=r_1(\xi_2)$\phantom{$r_3r_3r_3r_3r_3r_3$}\vphantom{$\int\limits_a$}}
\end{tabular}
&2\\
\hline

\end{tabular}
\end{center}

\end{table}


\section{General description of subgroups}
\label{General}

Given a Euclidean reflection group $G$, we present an algorithm to find
all reflection subgroups of $G$.

Let $F$ be a Coxeter simplex, and $P$ be a compact Coxeter polytope
in $\E^n$, where $G_{P}\subset G_{F}$.
Let  $P_1,...,P_s$ be indecomposable factors of $P$:
$P=P_1\times...\times P_s$, where $P_i$ are Euclidean simplices.
For each diagram $\Sigma(P_i)$, $i=1,...,s$,
take an arbitrary node $u_i$.
Then
$G_{P_1\setminus\mu(u_1)}\times...\times
G_{P_s\setminus\mu(u_s)}$ is a stabilizer $\mbox{Fix}(U,G_P)$
of some vertex $U$ of $P$ in  $G_{P}$.
Hence,  $\mbox{Fix}(U,G_F)$ contains $G_{P_1\setminus\mu(u_1)}\times...\times
G_{P_s\setminus\mu(u_s)}$.
Without loss of generality we may assume that $U$ is a vertex of $F$.

Therefore, any finite index reflection subgroup of $G_{F}$
can be obtained as a result of the following procedure:

1) Choose a vertex $U$ of $F$ and a maximal rank reflection subgroup
   $H$ of $\mbox{Fix}(U,G_F)$ (i.e. $U$ is the only
   point of $\E^n$ fixed by $H$).
   Denote by $K$ a fundamental cone of $H$.
   We have $K=K_1\times...\times K_s$, where $K_i$ are indecomposable
   cones.

2) For each cone $K_i$  take a mirror $\mu_i$ of  $G_F$
   such that $\mu_i\cap K_j =\emptyset$, $j\ne i$,
   and  $\mu_i$ cuts an acute-angled  simplex $P_i$ out of $K_i$.
   Let  $P=P_1\times...\times P_s$.
   Then $G_P$ is a finite index subgroup of $G_F$.

To give an explicit description of mirrors $\mu_i$ we use affine
root systems.

Let $\Delta_F$ be an affine root system such that
$G_F$ is the Weyl group of $\Delta_F$.
More precisely, for each multiple edge of $\Sigma(F)$ put an arrow
to arrange a Dynkin diagram $S(F)$.
We choose direction of arrows in order to obtain a diagram
contained in Table ``Aff1'' (see~\cite[Ch. 4]{Kac}).
The root system $\Delta_F$ consists of two disjoint parts: the set of
real roots  $\Delta^{\mbox{\small re}}_F$
and the set of  imaginary roots
$\Delta^{\mbox{\small im}}_F=\{\pm \delta,\pm 2\delta,...\}$ (see~\cite[Prop.
5.10]{Kac}).  Here $\delta=\sum a_i\alpha_i$, $\alpha_i$ are
simple roots of $\Delta_F$, and $a_i$ are the coefficients of the linear
dependency between the columns of generalized Cartan matrix.
Further, let $V$
be a special vertex of $F$, and $M$ be a set of mirrors of $G_F$
containing $V$.  Let $\Delta'_{F}\subset \Delta_{F}$ be a finite root
system that consists of all roots vanishing on mirrors contained in $M$.
By Prop.~6.3.  of~\cite{Kac}, $\Delta^{\mbox{\small re}}_F=\{ \alpha+n\delta \, |\,
\alpha \in \Delta'_{F}, n\in \Z \} $.

Let $P$ be a compact Coxeter polytope in $\E^n$, and
$G_P\subset G_F$. Consider roots of $\Delta_{F}$ vanishing on the
facets of $P$. These roots compose a set of simple roots for some root system
$\Delta \subset \Delta_{F}$. The Weyl group of $\Delta$ coincides with $G_P$.

Following the procedure described above, consider an arbitrary maximal rank
finite root system $\Delta'\subset \Delta_{F}$.
Let  $\Delta'=\Delta'_1+...+\Delta'_s$,  where $\Delta'_i$,
$i=1,...,s$, are indecomposable  components, and let $\Pi_i$ be a 
set of simple roots of $\Delta'_i$.
For each of $\Pi_i$ we should add a root
$\beta_i$ vanishing on $\mu_i$ and satisfying the foolowing two
conditions: $\beta_i$ is orthogonal to each root of $\Delta'_j$ if
$i\ne j$; for any $\gamma\in \Pi_i$ the angle formed up by
$\beta_i$ and  $\gamma$ is not acute.

  We can always take  $\theta_i+k_i\delta$ as $\beta_i$,
where $\theta_i$ is the lowest root of $\Delta'_i$, and $k_i\in \N$.
However, sometimes there exist additional roots satisfying the conditions
above.
In more details, let $S_i$ be a Dynkin diagram of $\Pi_i$.
If $S_i=B_l, C_l\ (l\ge 3), F_4, G_2$ or $C_2$,
there  exists a family of mutually parallel mirrors
such that any of these mirrors cuts an acute-angled polytope out of
$K_i$. Namely, the roots vanishing on these family equal
$\theta'+k\delta$, where   $k\in \N$, and $\theta'$ are listed
 in Table~\ref{add}.

\begin{table}
\tcaption{An additional root $\theta'$.}
\label{add}
\begin{center}
\begin{tabular}{|cc@{\quad}|c|}
\hline
\multicolumn{2}{|c|}{$S_i$\vphantom{$\int\limits_a^A$}}&$\theta'$\\
\hline
\raisebox{8pt}[10pt][9pt]{$B_l$}&
\raisebox{-0pt}[30pt][9pt]{\psfrag{2}{\scriptsize $1$}
\psfrag{3}{\scriptsize $2$}
\psfrag{4}{\scriptsize $l-1$}
\psfrag{5}{\scriptsize $l$}
\epsfig{file=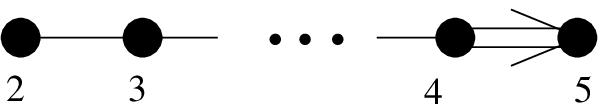,width=0.31\linewidth}}
&
\raisebox{9pt}[10pt][9pt]{$\frac{1}{2}(-\gamma_1+\theta)$}
\\
\hline
\raisebox{-2pt}[20pt][12pt]{$C_l$}&
\raisebox{-7pt}[10pt][12pt]{\psfrag{2}{\scriptsize $1$}
\psfrag{3}{\scriptsize $2$}
\psfrag{4}{\scriptsize $l-1$}
\psfrag{5}{\scriptsize $l$}
\epsfig{file=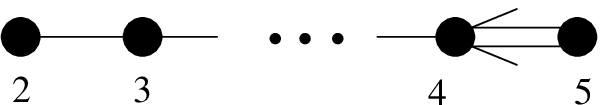,width=0.31\linewidth}}
&
\raisebox{-1pt}[10pt][12pt]{$\gamma_1+\theta$}
\\
\hline
\raisebox{-3pt}[15pt][12pt]{$G_2$}&
\raisebox{-7pt}[10pt][12pt]{\psfrag{3}{\scriptsize $1$}
\psfrag{4}{\scriptsize $1$}
\psfrag{5}{\scriptsize $2$}
\epsfig{file=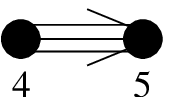,width=0.082\linewidth}}
&\raisebox{-2pt}[10pt][12pt]{$\gamma_1+\gamma_2+\theta$}
\\
\hline
\raisebox{-3pt}[15pt][12pt]{$C_2$}&
\raisebox{-7pt}[10pt][12pt]{\psfrag{3}{\scriptsize $1$}
\psfrag{4}{\scriptsize $1$}
\psfrag{5}{\scriptsize $2$}
\epsfig{file=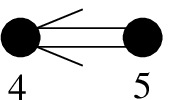,width=0.082\linewidth}}
&\raisebox{-2pt}[10pt][12pt]{$\gamma_1+\theta$}
\\
\hline
\raisebox{-3pt}[15pt][12pt]{$F_4$}&
\raisebox{-7pt}[10pt][12pt]{\psfrag{3}{\scriptsize $1$}
\psfrag{4}{\scriptsize $2$}
\psfrag{5}{\scriptsize $3$}
\psfrag{6}{\scriptsize $4$}
\epsfig{file=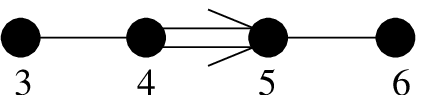,width=0.2\linewidth}}
&\raisebox{-2pt}[10pt][12pt]{$\gamma_1+\gamma_2+\gamma_3+\theta$}
\\
\hline
\end{tabular}
\end{center}
\end{table}

It is easy to check that there are no other possibilities for
 $\beta_i$.

\vspace{8pt}

Given a pair of Coxeter polytopes $F$ and $P$ one can ask
if it is possible that $G_P\subset G_F$.
The following theorem gives a criterion in terms of
Coxeter diagrams $\Sigma(F)$ and $\Sigma(P)$.

\begin{theorem}\label{general}
Let $F$ and  $P$ be Euclidean Coxeter polytopes.
A polytope $T$ satisfying $\Sigma(P)=\Sigma(T)$ and $G_{T}\subset G_F$
exists if and only if there exists an embedding of $G'_{P}$ into
$G'_{F}$.

\end{theorem}

\begin{proof}

To prove that the condition is necessary, assume that $T=P$. Since
$G_{P}\subset G_{F}$,  we may also assume that $P$ contains
$F$, and a special vertex $V$ of $P$ is also a vertex~of~$F$.
Then $G'_{P}(V)$ is a subgroup of the stabilizer $\mbox{Fix}(V,G_{F})$
of $V$ in $G_F$.
Clearly, there exists an embedding of
$\mbox{Fix}(V,G_{F})$ into  $G'_{F}$ (that takes any mirror $\mu\in
\mbox{Fix}(V,G_{F})$ to
a mirror parallel to $\mu$ and containing some fixed special vertex of
$F$).  Hence, we found an embedding $G'_{P}\hookrightarrow G'_{F}$.

To prove that the condition is sufficient, consider an image $G$
of $G'_P$ under the embedding $G'_{P}\hookrightarrow G'_{F}$.
Let $K$ be a fundamental chamber of $G$,
and $\Delta'\subset \Delta_{F}$ be a finite root system
whose simple roots vanish on the facets of $K$.
Let $\Delta'=\Delta'_1+...+\Delta'_s$,  where $\Delta'_i$,
$i=1,...,s$, are indecomposable components, and
let $\Pi_i$ be a set of simple roots of $\Delta'_i$.
For each of $\Pi_i$ we add a root
$\theta_i+\delta$, where $\theta_i$ is the lowest root of $\Delta'_{i}$.
Denote by $T_1$ a fundamental chamber of the Weyl group of the
direct sum of resulting root systems.  Clearly, $G_{T_1}\subset G_F$.

By the construction of $T_1$, the diagram $\Sigma(T_1)$ is
very similar to $\Sigma(P)$. The only possible difference
is that  some indecomposable components  $\widetilde B_l$ may be
substituted by $\widetilde C_l$ (and some  $\widetilde C_l$  may be
substituted by  $\widetilde B_l$).
For each of these components we take an index 2
(index  $2^{l-1}$ respectively) reflection subgroup described in
Table~\ref{min}.  In this way we obtain a subgroup with fundamental
chamber whose Coxeter diagram coincides with $\Sigma(P)$, and the
theorem is proved.

\end{proof}

\subsection*{Remark on the uniqueness}

In general, a decomposable subgroup $G_P\subset G_F$ is not determined by
$(\Sigma(P),\Sigma(F))$ and $[G_F:G_P]$.
For example, see Fig.~\ref{example} for fundamental polygones of three
different subgroups $G_P\subset G_F$ with
$\Sigma(F)=\widetilde C_2, \Sigma(P)=2\widetilde A_1$, $[G_F:G_P]=2\cdot
4$. Clearly, none of these subgroups is equivalent to another modulo the
automorphism group of $G_F$.
However, some additional conditions imply
the uniqueness (up to automorphism of $G_F$).

\begin{figure}[htb!]
\begin{center}
\epsfig{file=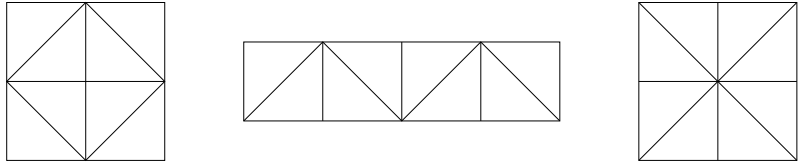,width=0.7\linewidth}\\
\end{center}
\caption{Fundamental domains of three subgroups with
$\Sigma(F)=\widetilde C_2, \Sigma(P)=2\widetilde A_1$ and
$[G_F:G_P]=2\cdot 4$. } \label{example} \end{figure}

Let $G_P\subset G_F$ be a subgroup, and let $M$ be a set of linear parts
of all elements of $G_P$. Suppose that $G_P$
contains reflections in all those mirrors of $G_F$ whose linear
part is contained in $M$.
Then we say that $G_P$ is {\it a block-maximal} subgroup of $G_F$.

\begin{lemma}
\label{block}
A block-maximal subgroup $G_P\subset G_F$ is determined 
by $(\Sigma(P),\Sigma(F))$ and $[G_F:G_P]$ uniquely 
up to an automorphism of $G_F$.

\end{lemma}

\begin{proof}

Without loss of generality we may assume that $F$ is contained in $P$.
Let $V$ be a special vertex of $F$.
For any mirror $\mu$ of $G_P$ there exists a mirror
containing $V$ and parallel to $\mu$.
Since the subgroup $G_P\subset G_F$ is block-maximal,
any of these mirrors is contained in $G_P$, so $V$
is a special  vertex of $P$.

By the definition of a block-maximal subgroup,
$G_P$ is completely determined by the stabilizer $G'_P$ of the special
vertex $V$. As it is shown in Section~\ref{sph}, an embedding of
the maximal rank finite subgroup $G'_P$ into $G'_F$ is usually unique
(up to an automorphism of $G'_F$). The only exclusions are the subgroups
described in Cor.~\ref{b_n-full}. However, for these subgroups the indices
$[G_F:G_P]$ differ by factor $2^k$, where $k\ne 0$ is the number defined
in Cor.~\ref{b_n-full}.

Suppose that $\Sigma(F)\ne \widetilde C_2,\widetilde G_2$ and
$\widetilde F_4$.  Then any automorphism of $G'_F$ can be extended to an
automorphism of $G_F$.  Hence, the subgroup $G_P\subset G_F$ is determined
 by $(\Sigma(P),\Sigma(F))$ uniquely up to an automorphism of $G_F$.

Now, suppose that  $\Sigma(F)= \widetilde C_2,\widetilde G_2$ or
$\widetilde F_4$. For each of these cases there exists the only
symmetry of $\Sigma'(F)$ that can not be extended to a symmetry of
$\Sigma(F)$. In other words, there are only two different embeddings of
$G'_P$  into $G'_F $ (up to an automorphism of $G_F$).
It is easy to check that for these embeddings the indices $[G_F:G_P]$
differ by factor  2, 3 and 4  respectively for
$\Sigma(F)= \widetilde C_2,\widetilde G_2$ and
$\widetilde F_4$.

\end{proof}

\section{Decomposable maximal subgroups}
\label{Max}

In this section we describe decomposable maximal reflection subgroups.

Let $G$ be a reflection group, $H\subset G$
be a finite reflection subgroup of $G$, and $M$ be a set of all mirrors of
$G$ that are parallel to mirrors of $H$.  Let $G_M$ be a group generated
by all reflections with respect to the mirrors contained in $M$.  Then we
say that $G_M$ is {\it a $G$-extension } of $H$. Obviously, $G_M$ is
a block-maximal subgroup of $G$.

\begin{lemma}
\label{general-max}

Let $F$ be a Coxeter simplex in  $\E^n$, $P$ be a decomposable polytope
in  $\E^n$, and  $G_{P}\subset G_F$ be a maximal reflection subgroup.
Then $G_P$ is a $G_F$-extension of some
 maximal reflection subgroup of $G'_F$.

Conversly, a $G_F$-extension of decomposable maximal reflection
subgroup of $G'_F$ is a maximal reflection subgroup of $G_F$.

\end{lemma}

\begin{proof}
Let $L$ be the set of all mirrors of
$G_F$ that are parallel to mirrors of $G_P$.
Denote by $G_L$ the group generated by
reflections in all mirrors contained in $L$.  Since $G_{P}\subseteq
G_{L}\subseteq G_F$ and the subgroup $G_{P}\subset G_F$ is maximal, either
$G_L=G_P$ or $G_L=G_F$.  The case  $G_L=G_F$ is impossible, since $G_P$ is
decomposable and $G_F$ is not. Thus, $G_L=G_P$ and $G_P$ is a
 $G_F$-extension of $G'_P$.

Furthermore, a spherical subgroup $G'_P$ is a maximal finite subgroup
of  $G'_F$ (indeed, if there exists a reflection subgroup $H$
satisfying $G'_P\subset H \subset G'_F$, then the $G_F$-extension of $H$
is a subgroup of $G_F$ containing $G_P$).
Therefore, any decomposable maximal subgroup is a $G_F$-extension of
some maximal finite (decomposable) subgroup of $G_F$.

From the other hand, it is clear that
a $G_F$-extension of any maximal finite subgroup
is a maximal subgroup of $G_F$.

\end{proof}

\begin{col}
A decomposable maximal reflection subgroup $G_P$ of an
 indecomposable reflection group  $G_F$ is block-maximal.
 Such a subgroup $G_P\subset G_F$  is determined  by
diagrams $\Sigma(F),\,\Sigma(P)$ and the index  $[G_F:G_P]$ uniquely 
up to an automorphism of  $G_F$.

\end{col}

\begin{proof}
By Lemma~\ref{general-max}, the subgroup $G_P\subset G_F$ is a
$G_F$-extension of some maximal subgroup. Thus, $G_P$ is block-maximal.
The uniquness follows from Lemma~\ref{block}.

\end{proof}

Using Table~12 of~\cite{Dyn} (see also~\cite{On}) and results of
Section~\ref{sph}, it is easy to find all maximal
rank decomposable maximal reflection subgroups of indecomposable finite
reflection groups.  Table~\ref{index} contains the complete list of them.
The same table contains the list of maximal decomposable Euclidean
reflection subgroups of indecomposable Euclidean reflection groups.  To
find the indices $[G_F:G_P]$ of these subgroups we calculated volumes of
$P$ and $F$:  $[G_F:G_P]=\frac{\mbox{Vol}(P)}{\mbox{Vol}(F)}$, 
where $\mbox{Vol}(T)$ is the
volume of $T$.

\begin{table}
\tcaption{Maximal rank decomposable maximal reflection subgroups
of spherical and Euclidean indecomposable reflection groups.
\vphantom{$\int\limits_a$}}
\label{index}
\begin{center}
\begin{tabular}{|c|c|c|c|c|c|}
\hline
\vphantom{$\int\limits^a_a$}$\Sigma '(F)$&$\Sigma '(P)$&$[G'_{F}:G'_{P}]$&
$\Sigma(F) $&$\Sigma(P) $&$[G_{F }:G_{P }]$\\
\hline
\vphantom{$\int\limits^a$}\raisebox{-10pt}[10pt][12pt]{$B_n=C_n
$}&\raisebox{-3pt}[10pt][12pt]{$B_k+B_{n-k}
$}&\raisebox{-10pt}[10pt][12pt]{${n\choose k}$}\vphantom{$\int\limits^A$}&
$\widetilde B_n$&$\widetilde B_k+\widetilde B_{n-k}$&$2{n\choose k}$\\
&\raisebox{4pt}[10pt][12pt]{($=C_k+C_{n-k}$)}&\vphantom{$\int\limits$}&
$\widetilde C_n$&$\widetilde C_k+\widetilde C_{n-k}$&$ {n\choose k}$\\
$D_n$&$D_k+D_{n-k}$&$2{n\choose k}$\vphantom{$\int\limits$}&
$\widetilde D_n$&$\widetilde D_k+\widetilde D_{n-k}$&$2^2{n\choose k}$
 \\
$G_2$&$2A_1$&$3$\vphantom{$\int\limits^a$}&
$\widetilde G_2$&$2\widetilde A_1$&$2\cdot 3$\\
$F_4$&$2A_2$&$2^5$\vphantom{$\int\limits^A$}&
$\widetilde F_4$&$2\widetilde A_2$&$2^5\cdot 3$\\
$E_6$&$A_5+A_1$&$2^2\cdot 3^2$\vphantom{$\int\limits^A$}&
$\widetilde E_6$&$\widetilde A_5+\widetilde A_1$&$2^3\cdot 3^2$\\
$E_6$&$3A_2$&$2^4\cdot 5$\vphantom{$\int\limits^A$}&
$\widetilde E_6$&$3\widetilde A_2$&$2^4\cdot 3^2 \cdot 5$\\
$E_7$&$D_6+A_1$&$3^2\cdot 7$\vphantom{$\int\limits^A$}&
$\widetilde E_7$&$\widetilde D_6+\widetilde A_1$&$2\cdot 3^2\cdot 7$\\
$E_7$&$A_5+A_2$&$2^5\cdot 3\cdot 7$\vphantom{$\int\limits^A$}&
$\widetilde E_7$&$\widetilde A_5+\widetilde A_2$&$2^5\cdot 3^2\cdot 7$\\
$E_8$&$E_7+A_1$&$2^3\cdot 3\cdot 5$\vphantom{$\int\limits^A$}&
$\widetilde E_8$&$\widetilde E_7+\widetilde A_1$&$2^4\cdot 3\cdot 5$\\
$E_8$&$E_6+A_2$&$2^6\cdot 5\cdot 7$\vphantom{$\int\limits^A$}&
$\widetilde E_8$&$\widetilde E_6+\widetilde A_2$&$2^6\cdot 3\cdot 5\cdot 7$\\
$E_8$&$2A_4$&$2^8\cdot 3^3\cdot 7$\vphantom{$\int\limits_A^A$}&
$\widetilde E_8$&$2\widetilde A_4$&$2^8\cdot 3^3\cdot 5\cdot 7$\\
\hline
\end{tabular}
\end{center}
\end{table}

\vspace{7pt}

\noindent
{\bf Remark.}
Table~\ref{index} shows  that $[G_{F }:G_{P }]$ is a multiple of
$[G'_{F}:G'_{P}]$. The reason of this is the following.
Let $V$ be a special vertex of $F$
and $O_V$ be the orbit of $V$ under the action of $G_F$.
Let $V_1,...,V_k$ be the points of $O_V$ contained either in $P$ or at the
boundary of $P$.
Let $G_{P(i)}$ be  the group generated by reflections in those facets of
$P$ that contain $V_i$.  To find the index $[G_{F }:G_{P }]$ it is
sufficient to calculate the number of images of $F$ under $G_F$ contained
in the polytope $P$:  $$ [G_{F }:G_{P
}]=\sum\limits_{i=1}^{k}\frac{|G'_{F}|}{|G_{P(i)}| }, $$ where $|G|$ is
the order of $G$.  Since $G_{P(i)}$ is a subgroup of $G'_{P}$, we have
$$|G_{P(i)}|=\frac{|G'_{P}|}{[G'_{P}:G_{P(i)}]}.$$
Hence,  $$\frac{[G_{F }:G_{P }]}{[G'_{F }:G'_{P}]}
=\sum\limits_{i=1}^{k} [G'_{P}:G_{P(i)}]\in \Z.$$

\section{Infinite index subgroups}
\label{Inf}

In previous sections we assumed that the polytope $P$ is compact,
so $G_P\subset G_F$ is a finite index subgroup.
However, sometimes a fundamental chamber of a discrete group generated
by reflections is not compact.
In this section, we discuss infinite index reflection subgroups
of discrete Euclidean indecomposable reflection group $G_F$ (where $F$ is a
Coxeter simplex).

Let $W$ be a discrete  group generated by reflections in $\E^n$, and $P$
be a fundamental chamber of $W$. Then $P$ is a generalized Coxeter
polytope, which is a convex domain bounded by finite number of hyperplanes
$f_1,...,f_k$, where  either $f_i$ is  parallel to $f_j$,
or  $f_i$ and $f_j$ form up an angle  $\frac{\pi}{m_{ij}}$, $m_{ij}\in\N$.
As it is shown in~\cite{Cox}, a generalized Euclidean Coxeter polytope
is a direct product of several simplices and simplicial cones.
A Coxeter diagram of this polytope is a union of several connected
parabolic and elliptic diagrams.

Let
$P=P_1\times...\times P_s\times P_{s+1}\times...\times P_{s+t}$
be a decomposition of a Coxeter polytope $P$ into indecomposable
components, where $P_1,...,P_s$ are Euclidean simplices
and $P_{s+1},...,P_{s+t}$ are indecomposable simplicial cones.
Then $\Sigma(P)$ is a union of $s$ connected parabolic diagrams  and $t$
connected elliptic diagrams.
Let $v_1,...,v_s$ be special vertices of the diagrams
$\Sigma(P_1),...,\Sigma(P_s)$ respectively.

\vspace{8pt}
A direct generalization of arguments of
Section~\ref{General} leads to the following description of subgroups of
given indecomposable Euclidean reflection group.

Consider a simplex $F$ in $\E^n$ and a root system $\Delta_F$
described in Section~\ref{General}.
Then any reflection subgroup (that may be of infinite index) of $G_F$
can be obtained as a result of the following procedure:
choose a finite reflection subgroup  $G\subset G_F$
($G$ may not be of maximal rank);
for each indecomposable component  $G_i$ of $G$ consider the
corresponding root system $\Delta_i\subset \Delta_{F}$
and take a positive integer $k_i$.
Now enlarge some of $\Delta_i$ by the roots $\theta_i+k_i\delta$,
where $\theta_i$ is the lowest root of $\Delta_i$
(one can take $\theta'$ instead of $\theta$ for
some components, see Table~\ref{add}),
the rest systems $\Delta_i$ leave unchanged.
The Weyl group $W$ of the resulting root system is a reflection
subgroup of $G_F$  ($W$ may be of infinite index).

\vspace{10pt}

Now, suppose that the subgroup $W\subset G_F$ is maximal.
Finite index maximal reflection subgroups are classified in
Sections~\ref{simpl} and~\ref{Max}.
In the following theorem we list all infinite index maximal reflection
subgroups.

\begin{theorem}
Let $F$ be a Euclidean Coxeter simplex
and $W\subset G_F$ be an infinite index  maximal reflection subgroup.
Let $\Sigma$ be the Coxeter diagram of fundamental chamber of $W$.
Then $(\Sigma(F),\Sigma)$ coincides with one of the following pairs:
$(\widetilde A_n,\widetilde A_k+\widetilde A_{n-1-k})$, $k=0,1,2,...,n-1$,
$(\widetilde D_n,\widetilde D_{n-1})$,
$(\widetilde D_n,\widetilde A_{n-1})$,
$(\widetilde E_6,\widetilde D_5)$ and
$(\widetilde E_7,\widetilde E_6)$.
\end{theorem}

\begin{proof}
Let $W'\subset G'_F$ be a finite maximal reflection subgroup of
$W$, and  $\Sigma'$ be a Coxeter diagram of fundamental domain of $W'$.
Since $W\subset G_F$ is a maximal subgroup,
the subgroup $W'\subset G'_F$ is also maximal.
Results of Section~\ref{sph} imply that the groups $B_n$, $F_4$ and
$G_2$ have no maximal subgroups of non-maximal rank.
By Cor.~\ref{sl}, all maximal reflection subgroups of non-maximal rank
are listed in~Table~12 of~\cite{Dyn}. Namely,
in this case the pair $(\Sigma'(F),\Sigma')$ coincides with
one of the following:  $(A_n,A_k+A_{n-1-k})$, $k=0,1,2,...,n-1$,
$(D_n,D_{n-1})$, $(D_n,A_{n-1})$, $(E_6,D_5)$ and $(E_7,E_6)$.  This
proves the theorem.

\end{proof}


\vspace{20pt}

\noindent
{\it Independent University of Moscow\\
 E-mail addresses:\quad felikson@mccme.ru\\
\phantom{E-mail addresses:}\quad pasha@mccme.ru
}

}
\end{document}